\documentclass{amsart}
\usepackage{amsmath,amssymb, amsbsy}
\usepackage[dvips]{graphicx}
\newcommand{\R}{{\mathbb R}}
\newcommand{\N}{{\mathbb N}}
\newcommand{\Z}{{\mathbb Z}}

\newcommand{\p}{\partial}

\newcommand{\iy}{\infty}

\newcommand{\dyle}{\displaystyle}
\newcommand{\dint}{\dyle\int}
\newcommand{\weakly}{\rightharpoonup}
\newcommand{\e }{\varepsilon}

\renewcommand{\ge }{\geqslant}
\renewcommand{\geq }{\geqslant}

\renewcommand{\leq }{\leqslant}

\newenvironment{pf}{\noindent{\sc Proof}.\enspace}{\hfill\qed\medskip}
\newenvironment{pfn}[1]{\noindent{\bf Proof of
    {#1}.\enspace}}{\hfill\qed\medskip}
\newtheorem{Theorem}{Theorem}[section]
\newtheorem{Corollary}[Theorem]{Corollary}
\newtheorem{Lemma}[Theorem]{Lemma}
\newtheorem{Proposition}[Theorem]{Proposition}
\theoremstyle{definition} \newtheorem{Definition}[Theorem]{Definition}

\newtheorem{remark}[Theorem]{Remark}
\begin{document}

\title[Heat equations with
inverse-square potentials]
{Classification of  local asymptotics
for solutions to   heat equations with
inverse-square potentials}

\author[Veronica Felli \and Ana Primo]{Veronica Felli \and Ana Primo}
\address{\hbox{\parbox{5.7in}{\medskip\noindent{Universit\`a di Milano
        Bicocca,\\
        Dipartimento di Ma\-t\-ema\-ti\-ca e Applicazioni, \\
        Via Cozzi
        53, 20125 Milano, Italy. \\[3pt]
        \em{E-mail addresses: }{\tt veronica.felli@unimib.it,
         ana.primo@unimib.it}.}}}}

\date{\today}

\thanks{
  The second author is supported by project MTM2007-65018, MEC, Spain}

 \begin{abstract}
   \noindent Asymptotic behavior of solutions to heat equations with spatially
   singular inverse-square potentials is studied. By combining a
   parabolic Almgren type monotonicity formula with blow-up methods, we
   evaluate the exact behavior near the singularity of solutions to
   linear and subcritical semilinear parabolic equations with  Hardy type
   potentials.  As a remarkable byproduct, a unique continuation
   property is obtained.
 \end{abstract}

\maketitle

\section{Introduction and statement of the main results}\label{intro}

We aim to describe the asymptotic behavior
near the singularity of solutions to backward evolution equations
with inverse square singular potentials
of the form
\begin{equation}\label{prob}
u_t+\Delta u+\dfrac{a(x/|x|)}{|x|^2}\,u+f(x,t, u(x,t))=0,
\end{equation}
in $\R^N\times (0,T)$, where $T>0$, $N\geq 3$, $a\in
L^{\infty}({\mathbb S}^{N-1})$ and $f:\R^N\times(0,T)\times\R\to \R$.
Inverse square potentials are related to the well-known classical
Hardy's inequality
\begin{equation*}
  \int_{\R^N}\,|\nabla u(x)|^{2}\,dx\geq
\bigg(\frac{N-2}{2}\bigg)^{\!\!2}\int_{\R^N}
  \frac{u^{2}(x)}{|x|^{2}}\,dx,
  \quad\text{for all }u\in\mathcal{C}_0^\infty(\R^N),\quad N\ge 3,
\end{equation*}
see e.g. \cite{GP,HLP}.
Parabolic problems with singular inverse square Hardy potentials arise
in the linearization of standard combustion models, see \cite{PV}.
The properties of the heat operator are strongly affected by the
presence of the singular inverse square potential, which, having the
same order of homogeneity as the laplacian and failing to belong to
the Kato class, cannot be regarded as a lower order term. Hence,
singular problems with inverse square potentials represent a
borderline case with respect to the classical theory of parabolic
equations. Such a criticality makes parabolic equations of type
(\ref{prob}) and their elliptic versions quite challenging from the
mathematical point of view, thus motivating a large literature which,
starting from the pioneering paper by \cite{BaGo}, has been devoted to
their analysis, see e.g. \cite{GP,vazquez_zuazua} for the parabolic
case and \cite{AFP,smets,terracini} for the elliptic counterpart.  In
particular, the influence of the Hardy potential in semilinear
parabolic problems has been studied in \cite{APP}, in the case
$f(x,t,s)=s^p$, $p>1$, and for $a(x/|x|)=\lambda$, $\lambda>0$; the
analysis carried out in \cite{APP} highlighted a deep difference with
respect to the classical heat equation ($\lambda=0$), showing that, if
$\lambda>0$, there exists a critical exponent $p_+(\lambda)$ such that for
$p\ge p_+(\lambda)$, there is no solution even in the weakest sense for any
nontrivial initial datum.

The present paper is addressed to the problem of describing the
behavior of solutions along the directions $(\lambda x,\lambda^2 t)$
naturally related to the heat operator. Indeed, the unperturbed
operator $u_t+\Delta u+\frac{a(x/|x|)}{|x|^2}\,u$ is invariant under
the action $(x,t)\mapsto(\lambda x,\lambda^2 t)$. Then we are
interested in evaluating the asymptotics of
$$
u(\sqrt{t} x,t)\quad\text{as }t\to 0^+
$$
for solutions to (\ref{prob}). Our analysis will show that $u(\sqrt{t}
x,t)$ behaves as a singular self-similar eigenfunction of the
Ornstein-Uhlenbeck operator with  inverse square potential,
multiplied by a power of $t$ related to the corresponding eigenvalue,
which can be selected by the limit of a frequency type function
associated to the problem.

We consider both linear and subcritical
semilinear parabolic equations of type (\ref{prob}). More precisely,
we deal with the case $f(x,t,s)= h(x,t)s$ corresponding to the linear
problem
\begin{equation}\label{prob1}
  u_t+\Delta u+\dfrac{a(x/|x|)}{|x|^2}\,u+h(x,t)u=0,\quad
\text{in }\R^N\times (0,T),
\end{equation}
with a perturbing potential $h$ satisfying
\begin{equation}\label{eq:der}
h,h_t\in L^{r}\big((0,T),L^{{N}/{2}}(\R^N)\big)
\quad\text{for some }r>1,
\quad h_t\in L^{\infty}_{\rm loc}\big((0,T),L^{{N}/{2}}(\R^N)\big),
\end{equation}
and negligible with respect to the inverse square potential $|x|^{-2}$
near the singularity in the sense that there exists $C_h>0$ such that
\begin{equation}\label{eq:h}
  |h(x,t)|\leq C_h(1+|x|^{-2+\e}) \quad
  \text{for all }t\in (0,T),\text{ a.e. }x\in\R^N,
  \text{ and for some }\e\in(0,2).
\end{equation}
We also treat the semilinear case  $f(x,t,s)= \varphi (x, t,
s)$, with a nonlinearity $\varphi\in C^1(\R^N\times(0,T)\times\R)$
 satisfying the following growth condition
\begin{equation}\label{eq:fi}
\begin{cases}
  \dfrac{|\varphi (x,t,s)|+|x\cdot\nabla_x \varphi (x,t,s)|
    +|t\frac{\p \varphi}{\p t}(x,t,s)|}{|s|}\leq C_\varphi (1+|s|^{p-1})\\[10pt]
\big|\varphi (x,t,s)-s{\textstyle{\frac{\p \varphi}{\p s}}}(x,t,s)\big|
\leq C_\varphi|s|^q
\end{cases}
\end{equation}
for all $(x,t,s)\in \R^N\times(0,T)\times\R$ and some $1< p<2^*-1$
and $2\leq q<p+1$,
where $2^*=\frac{2N}{N-2}$ is the critical exponent for Sobolev's
embedding and $C_\varphi>0$ is independent of $x\in\R^N$, $t\in (0,T)$,
and $s\in\R$. In particular, we are going to classify the behavior of
solutions to the semilinear parabolic problem
\begin{equation}\label{prob2}
u_t+\Delta u+\dfrac{a(x/|x|)}{|x|^2}\,u+\varphi (x, t, u(x,t))=0,
\quad
\text{in }\R^N\times (0,T),
\end{equation}
satisfying
\begin{equation}\label{eq:u1}
 u\in L^{\infty}(0,T, L^{p+1}(\R^N))
\end{equation}
and
\begin{equation}\label{eq:u2}
t  u_t\in L^{\infty}(0,T, L^{\frac{p+1}{p+1-q}}(\R^N)) \mbox { and }
  \sup\limits_{t\in (0,T)}t^{N/2}\int_{\R^N}
|x|^{\frac{2(p+1)}{p-1}}
   |u(\sqrt t x,t)|^{p+1}\,dx<\infty.
\end{equation}

In order to introduce a suitable notion of solution to
(\ref{prob}), for every
$t>0$ let us define the space ${\mathcal H}_t$ as
the completion of $C^{\infty}_{\rm c}(\R^N)$ with respect to
$$
\|u\|_{{\mathcal H}_t}=\bigg(\int_{\R^N}\big(t|\nabla
u(x)|^2+|u(x)|^2\big) G(x,t)\,dx\bigg)^{\!\!1/2},
$$
where
\begin{equation*}
G(x,t)=t^{-N/2}\exp\Big(-\frac{|x|^2}{4t}\Big)
\end{equation*}
is the heat kernel satisfying
\begin{equation}\label{eq:heatker}
 G_t-\Delta G=0\quad\text{and}\quad \nabla G(x,t)=-\frac{x}{2t}\,G(x,t)
\quad\text{in }\R^N\times (0,+\infty).
\end{equation}
We denote as $\big({\mathcal H}_t\big)^\star$ the dual space of
${\mathcal H}_t$ and by ${}_{({\mathcal H}_t)^\star}\langle
\cdot,\cdot\rangle_{{\mathcal H}_t}$ the corresponding duality product.
For every $t>0$, we also define  the space ${\mathcal L}_t$ as
the completion of $C^{\infty}_{\rm c}(\R^N)$ with respect to
$$
\|u\|_{{\mathcal L}_t}=\bigg(\int_{\\R^N}|u(x)|^2 G(x,t)\,dx\bigg)^{\!\!1/2}.
$$
\begin{Definition}\label{def:solution}
  We say that $u\in L^1_{\rm loc
  }(\R^N\times(0,T))$ is a weak solution to
  (\ref{prob}) in $\R^N\times(0,T)$ if
\begin{align}
  \label{eq:defsol1}&\int_\tau^T\|u(\cdot,t)\|^2_{{\mathcal
      H}_t}\,dt<+\infty,\quad\int_\tau^T\Big\|u_t+\frac{\nabla u\cdot
    x}{2t}\Big\|^2_{({\mathcal H}_t)^\star}dt<+\infty \text{ for all
  }\tau\in (0,T),\\
\label{eq:defsol2}&{\phantom{\bigg\langle}}_{{\mathcal
    H}_t^\star}\bigg\langle
u_t+\frac{\nabla u\cdot x}{2t},\phi
\bigg\rangle_{{\mathcal H}_t}\\
&\notag\qquad=
\int_{\R^N}\bigg(\nabla u(x,t)\cdot \nabla \phi(x)-
\dfrac{a(x/|x|)}{|x|^2}\,u(x,t)\phi(x)-f(x,t,u(x,t))\phi(x)\bigg)G(x,t)\,dx
\end{align}
for a.e. $t\in (0,T)$ and for each $\phi\in {\mathcal
    H}_t$.
\end{Definition}
It will be clear from the parabolic Hardy type inequality of Lemma
\ref{Hardytemp} and the Sobolev weighted inequality of Corollary
\ref{cor:ineqSob}, that the integral
$\int_{\R^N}f(x,t,u(x,t))\phi(x)G(x,t)dx$ in the above definition is
finite for a.e. $t\in(0,T)$, both in the linear case $f(x,t,s)=
h(x,t)s$ under assumptions (\ref{eq:der}--\ref{eq:h}) and in the
semilinear case $f(x,t,s)= \varphi (x, t, s)$ under condition
(\ref{eq:fi}) and for $u$ satisfying (\ref{eq:u1}).

\begin{remark}\label{rem:uv}
  If $u\in L^1_{\rm loc }(\R^N\times(0,T))$ satisfies (\ref{eq:defsol1}), then
the function
$$
v(x,t):=u(\sqrt{t}x,t)
$$
satisfies
\begin{equation}\label{eq:4}
v\in L^2(\tau,T;\mathcal H)\quad\text{and}\quad
v_t\in L^2(\tau,T;(\mathcal H)^\star) \quad\text{for all
  }\tau\in (0,T),
\end{equation}
where we have set
$$
\mathcal H:={\mathcal H}_1,
$$
i.e. ${\mathcal H}$ is
the completion of $C^{\infty}_{\rm c}(\R^N)$ with respect to
$$
\|v\|_{{\mathcal H}}=\bigg(\int_{\R^N}\big(|\nabla
v(x)|^2+|v(x)|^2\big) e^{-|x|^2/4}\,dx\bigg)^{\!\!1/2}.
$$
We notice that from (\ref{eq:4}) it follows that
$$
v\in C^0([\tau,T],{\mathcal L}),
$$
see e.g. \cite[Theorem 1.2]{SH}, where $\mathcal L:={\mathcal L}_1$ is
the completion of $C^{\infty}_{\rm c}(\R^N)$ with respect to the norm
$\|v\|_{{\mathcal L}}=\big(\int_{\R^N}|v(x)|^2e^{-|x|^2/4}\,dx\big)^{1/2}$.
Moreover the function
$$
t\in[\tau,T]\mapsto \|v(t)\|^2_{{\mathcal L}}=\int_{\R^N}u^2(x,t)G(x,t)\,dx
$$
is absolutely continuous  and
$$
\frac12\frac1{dt} \int_{\R^N}u^2(x,t)G(x,t)=
\frac12\frac1{dt}\|v(t)\|_{{\mathcal L}}^2=
{}_{{\mathcal H}^\star}\langle
v_t(\cdot,t),v(\cdot,t)
\rangle_{{\mathcal H}}
={\phantom{\bigg\langle}}_{{\mathcal
    H}_t^\star}\bigg\langle
u_t+\frac{\nabla u\cdot x}{2t},u(\cdot,t)
\bigg\rangle_{{\mathcal H}_t}
$$
for a.e. $t\in(0,T)$.
\end{remark}
\begin{remark}\label{rem:v2}
  If $u$ is a weak solution to (\ref{prob}) in the sense of definition
  \ref{def:solution}, then the function $v(x,t):=u(\sqrt{t}x,t)$
  defined in Remark \ref{rem:uv} is  a weak solution to
\begin{equation*}
  v_t+\frac1t\bigg(\Delta v-\frac x2\cdot \nabla v+
\dfrac{a(x/|x|)}{|x|^2}\,v+tf(\sqrt tx,t,v(x,t))\bigg)=0,
\end{equation*}
in the sense that, for every $\phi\in{\mathcal H}$,
\begin{multline}\label{eq:24}
{\phantom{\big\langle}}_{{\mathcal
    H}^\star}\!\big\langle v_t,\phi \big\rangle_{{\mathcal H}}\\=
\frac1t\int_{\R^N}\!\!\bigg(\!\nabla v(x,t)\!\cdot \!\nabla \phi(x)-
\dfrac{a\big(\frac{x}{|x|}\big)}{|x|^2}\,v(x,t)\phi(x) -t\,f(\sqrt
tx,t,v(x,t))\phi(x)\!\bigg)G(x,1)\,dx.
\end{multline}
In particular, if $u$ is a weak solution to (\ref{prob1}), then
 $v(x,t):=u(\sqrt{t}x,t)$ weakly solves
\begin{equation*}
  v_t+\frac1t\bigg(\Delta v-\frac x2\cdot \nabla v+
\dfrac{a(x/|x|)}{|x|^2}\,v+t h(\sqrt tx,t)v\bigg)=0,
\end{equation*}
whereas, if $u$ is a weak solution to (\ref{prob2}), then
$v(x,t):=u(\sqrt{t}x,t)$ weakly solves
\begin{equation*}
  v_t+\frac1t\bigg(\Delta v-\frac x2\cdot \nabla v+
\dfrac{a(x/|x|)}{|x|^2}\,v+t\varphi(\sqrt tx,t,v)\bigg)=0.
\end{equation*}
\end{remark}
We give a precise description of the asymptotic behavior at the singularity of
solutions to (\ref{prob1}) and (\ref{prob2}) in terms of
the eigenvalues and eigenfunctions of the Ornstein-Uhlenbeck operator
with singular inverse square potential
\begin{equation}\label{eq:13}
L:{\mathcal H}\to({\mathcal H})^\star,\quad
L=-\Delta+\frac{x}2\cdot\nabla-\frac{a(x/|x|)}{|x|^2},
\end{equation}
acting as
$$
{}_{{\mathcal H}^\star}\langle
Lv,w
\rangle_{{\mathcal H}}
=
\int_{\R^N}\bigg(\nabla
v(x)\cdot\nabla w(x)-\frac{a(x/|x|)}{|x|^2}\,v(x)w(x)\bigg)G(x,1)\,dx,
\quad\text{for all }v,w\in{\mathcal H}.
$$
In order to describe the spectrum of $L$, we consider the operator
$-\Delta_{\mathbb S^{N-1}}-a(\theta)$ on the unit $(N-1)$-dimensional
sphere $\mathbb S^{N-1}$. For any $a\in L^{\infty}\big({\mathbb
  S}^{N-1}\big)$, $-\Delta_{\mathbb S^{N-1}}-a(\theta)$ admits a
diverging sequence of eigenvalues
$$
\mu_1(a)<\mu_2(a)\leq\cdots\leq\mu_k(a)\leq \cdots,
$$
the first of which is simple and can be characterized as
\begin{equation}\label{eq:67}
\mu_1(a)=\min_{\psi\in H^1(\mathbb
S^{N-1})\setminus\{0\}}\frac{\int_{\mathbb S^{N-1}}|\nabla_{\mathbb
S^{N-1}}\psi(\theta)|^2\,dS(\theta)-\int_{\mathbb S^{N-1}}a(\theta)
\psi^2(\theta)\,dS(\theta)}{\int_{\mathbb S^{N-1}}\psi^2(\theta)\,dS(\theta)},
\end{equation}
see  \cite{FMT2}.
Moreover the quadratic form associated to $-\Delta
-\frac{a(x/|x|)}{|x|^2}$
 is positive definite if and only if
\begin{equation}\label{eq:posde}
\mu_1(a)>-\frac{(N-2)^2}4,
\end{equation}
see \cite[Lemma 2.5]{FMT2}.
  To each $k\in\N$, $k\geq 1$, we associate a
$L^{2}\big({\mathbb S}^{N-1}\big)$-normalized eigenfunction $\psi_k$
of the operator $-\Delta_{\mathbb S^{N-1}}-a(\theta)$ corresponding to the $k$-th
eigenvalue $\mu_{k}(a)$, i.e. satisfying
\begin{equation}\label{eq:2rad}
\begin{cases}
  -\Delta_{\mathbb S^{N-1}}\psi_k(\theta)-a(\theta)\psi_k(\theta)
  =\mu_k(a)\,\psi_k(\theta),&\text{in }{\mathbb S}^{N-1},\\[3pt]
  \int_{{\mathbb S}^{N-1}}|\psi_k(\theta)|^2\,dS(\theta)=1.
\end{cases}
\end{equation}
In the enumeration $\mu_1(a)<\mu_2(a)\leq\cdots\leq\mu_k(a)\leq \cdots$ we
repeat each eigenvalue as many times as its multiplicity; thus exactly
one eigenfunction $\psi_k$ corresponds to each index $k\in\N$.  We can
choose the functions $\psi_k$ in such a way that they form an
orthonormal basis of $L^2({\mathbb S}^{N-1})$.

The following proposition describes completely the spectrum of the
operator $L$, thus extending to the anisotropic case the spectral
analysis performed in \cite[\S 9.3]{vazquez_zuazua} in the isotropic
case $a(\theta)\equiv \lambda$; see also \cite[\S 4.2]{CH} and
\cite[\S 2]{ES} for the non singular case.

\begin{Proposition}\label{p:explicit_spectrum}
The set of the eigenvalues  of the operator $L$ is
$$
\big\{ \gamma_{m,k}: k,m\in\N, k\geq 1\big\}
$$
where
\begin{equation}\label{eq:65}
\gamma_{m,k}=m-\frac{\alpha_k}2, \quad
\alpha_k=\frac{N-2}{2}-\sqrt{\bigg(\frac{N-2}{2}\bigg)^{\!\!2}+\mu_k(a)},
\end{equation}
and $\mu_k(a)$ is the $k$-th eigenvalue of the operator
$-\Delta_{\mathbb S^{N-1}}-a(\theta)$ on the sphere $\mathbb
S^{N-1}$. Each eigenvalue $\gamma_{m,k}$ has finite multiplicity equal
to
$$
\#\bigg\{j\in\N,j\geq 1:
\gamma_{m,k}+\frac{\alpha_j}2\in\N\bigg\}
$$
and a basis of the corresponding eigenspace is
$$
\left\{V_{n,j}: j,n\in\N,j\geq
  1,\gamma_{m,k}=n-\frac{\alpha_j}2 \right\},
$$
where
\begin{equation}\label{eq:66}
V_{n,j}(x)=
|x|^{-\alpha_j}P_{j,n}\bigg(\frac{|x|^2}{4}\bigg)
\psi_j\Big(\frac{x}{|x|}\Big),
\end{equation}
$\psi_j$ is an eigenfunction of the operator
$-\Delta_{\mathbb S^{N-1}}-a(\theta)$ on the sphere $\mathbb S^{N-1}$
associated to the $j$-th eigenvalue $\mu_{j}(a)$ as in (\ref{eq:2rad}),
  and $P_{j,n}$
is the polynomial of degree $n$ given by
$$
P_{j,n}(t)=\sum_{i=0}^n
\frac{(-n)_i}{\big(\frac{N}2-\alpha_j\big)_i}\,\frac{t^i}{i!},
$$
denoting as $(s)_i$, for all $s\in\R$, the Pochhammer's symbol
$(s)_i=\prod_{j=0}^{i-1}(s+j)$, $(s)_0=1$.
\end{Proposition}
The following theorems provide a classification of singularity rating
of any solution $u$ to (\ref{prob}) based on the limit as $t\to 0^+$ of the
{\it Almgren type frequency function} (see \cite{almgren,poon}),
\begin{equation}\label{eq:64}
{\mathcal N}_{f,u}(t)=\frac{t
\int_{\R^N}\big(|\nabla u(x,t)|^2-
\frac{a(x/|x|)}{|x|^2}u^2(x,t)-f(x,t, u(x,t))u(x,t)\big)G(x,t)\,dx}
{\int_{\R^N}u^2(x,t)\, G(x,t)\,dx}.
\end{equation}
In the linear case $f(x,t,u)=h(x,t)u$, the behavior of weak solutions
to (\ref{prob1}) is described by the following theorem.
\begin{Theorem}\label{asym1}
  Let $u\not\equiv 0$ be a weak solution to (\ref{prob1}) in the sense
  of Definition \ref{def:solution}, with  $h$ satisfying
  (\ref{eq:der}) and (\ref{eq:h}) and $a\in L^{\infty}\big({\mathbb
    S}^{N-1}\big)$ satisfying (\ref{eq:posde}). Then there exist $m_0,k_0\in\N$,
  $k_0\geq 1$, such that
\begin{equation}\label{eq:671}
\lim_{t\to 0^+}{\mathcal N}_{hu, u}(t)=\gamma_{m_0,k_0},
\end{equation}
where ${\mathcal N}_{hu,u}$ is defined in (\ref{eq:64}) and
$\gamma_{m_0,k_0}$ is as in (\ref{eq:65}). Furthermore, denoting as
$J_0$ the finite set of indices
\begin{equation}\label{eq:79}
J_0=\{(m,k)\in\N\times(\N\setminus\{0\}):m-\frac{\alpha_{k}}2=\gamma_{m_0,k_0}\},
\end{equation}
for all $\tau \in(0,1)$ there holds
\begin{equation}\label{eq:80}
\lim_{\lambda\to0^+}\int_\tau^1
\bigg\|\lambda^{-2\gamma_{m_0,k_0}}u(\lambda x,\lambda^2t)
-t^{\gamma_{m_0,k_0}}\sum_{(m,k)\in J_0}\beta_{m,k}\widetilde V_{m,k}(x/\sqrt t)
\bigg\|_{{\mathcal H}_t}^2dt=0
\end{equation}
and
\begin{equation}\label{eq:81}
\lim_{\lambda\to0^+}\sup_{t\in[\tau,1]}
\bigg\|\lambda^{-2\gamma_{m_0,k_0}}u(\lambda x,\lambda^2t)
-t^{\gamma_{m_0,k_0}}\sum_{(m,k)\in J_0}\beta_{m,k}\widetilde V_{m,k}(x/\sqrt t)
\bigg\|_{{\mathcal L}_t}=0,
\end{equation}
where
$\widetilde V_{m,k}=
V_{m,k}/\|V_{m,k}\|_{\mathcal L}$,
$V_{m,k}$ are as in (\ref{eq:66}),
\begin{multline}\label{eq:beta1}
\beta_{m,k}=\Lambda^{-2\gamma_{m_0,k_0}}
\int_{\R^N}u(\Lambda x,\Lambda^2)\widetilde V_{m,k}(x)G(x,1)\,dx\\
+2\int_0^{\Lambda}s^{1-2\gamma_{m_0,k_0}}
\bigg(\int_{\R^N}h(s x, s^2) u(s x,s^2)\widetilde V_{m,k}(x)G(x,1)\,dx\bigg)
ds
\end{multline}
for all $\Lambda\in(0,\Lambda_0)$ and for some
$\Lambda_0\in(0,\sqrt{T})$, and $\beta_{m,k}\neq0$ for some $(m,k)\in J_0$.
\end{Theorem}

An analogous result holds in the semilinear case for solutions to
(\ref{prob2}) satisfying the further conditions (\ref{eq:u1}) and
(\ref{eq:u2}).
\begin{Theorem}\label{asym2}
  Let $a\in L^{\infty}\big({\mathbb S}^{N-1}\big)$ satisfy
  (\ref{eq:posde}) and $\varphi\in C^1(\R^N\times(0,T)\times\R)$ such
  that (\ref{eq:fi}) holds.  If $u\not\equiv 0$ satisfies
  (\ref{eq:u1}--\ref{eq:u2}) and is a weak solution to (\ref{prob2}) in the sense
  of Definition \ref{def:solution}, then there exist $m_0,k_0\in\N$,
  $k_0\geq 1$, such that
\begin{equation}\label{eq:672}
\lim_{t\to 0^+}{\mathcal N}_{\varphi,u}(t)=\gamma_{m_0,k_0},
\end{equation}
where ${\mathcal N}_{\varphi,u}$ is defined in (\ref{eq:64}) and $\gamma_{m_0,k_0}$ is as in
(\ref{eq:65}). Furthermore, letting $J_0$  the finite set of indices
defined in (\ref{eq:79}),
for all $\tau \in(0,1)$
convergences (\ref{eq:80}) and (\ref{eq:81}) hold with
\begin{multline}\label{eq:beta2}
\beta_{m,k}=\Lambda^{-2\gamma_{m_0,k_0}}
\int_{\R^N}u(\Lambda x,\Lambda^2)\widetilde V_{m,k}(x)G(x,1)\,dx\\
+2\int_0^{\Lambda}s^{1-2\gamma_{m_0,k_0}}
\bigg(\int_{\R^N}\varphi(s x, s^2, u(s x,s^2))\widetilde V_{m,k}(x)G(x,1)\,dx\bigg)
ds
\end{multline}
for all $\Lambda\in(0,\Lambda_0)$ and for some
$\Lambda_0\in(0,\sqrt{T})$, and $\beta_{m,k}\neq0$ for some $(m,k)\in J_0$.
\end{Theorem}
\eqref{eq:beta1} and \eqref{eq:beta2} can be seen as Cauchy's integral
type formulas for solutions to problems \eqref{prob1} and
\eqref{prob2}, since they allow reconstructing, up to the
perturbation, the solution at the singularity by the values it takes
at any positive time.

The proofs of theorems \ref{asym1} and \ref{asym2} are based on a
parabolic Almgren type monotonicity formula combined with blow-up
methods. Almgren type frequency functions associated to parabolic
equations were first introduced by C.-C. Poon in \cite{poon}, where
unique continuation properties are derived by proving a monotonicity
result which is the parabolic counterpart of the monotonicity formula
introduced by Almgren in \cite{almgren} and extended by Garofalo and
Lin in \cite{GL} to elliptic operators with variable coefficients.  A
further development in the use of Almgren monotonicity methods to study
regularity of solutions to parabolic problems is due to the recent
paper \cite{C}.
 We also mention that an
Almgren type monotonicity method combined with blow-up  was
used in \cite{FFT} in an elliptic context
to study the behavior of solutions to stationary Schr\"odinger
equations with singular electromagnetic potentials.

Theorem \ref{asym1} and Theorem \ref{asym2} imply {\it a strong unique
continuation property} at the singularity, as the following corollary
states.

\begin{Corollary}\label{cor:uniq_cont}
  Suppose that either $u$ is a weak solution to (\ref{prob1}) under
  the assumptions of Theorem \ref{asym1} or $u$ satisfies
  (\ref{eq:u1}--\ref{eq:u2}) and weakly solves (\ref{prob2}) under the
  assumptions of Theorem \ref{asym2}. If
\begin{equation}\label{eq:uniq_cont}
u(x,t)=O\big((|x|^2+t)^k\big)\quad\text{as }(x,t)\to(0,0)
\quad\text{for all }k\in\N,
\end{equation}
then $u\equiv 0$ in $\R^N\times(0,T)$.
\end{Corollary}

As a byproduct of the proof of Theorems \ref{asym1} and \ref{asym2},
we also obtain the following result, which can be regarded as a {\it unique
continuation property} with respect to time.

\begin{Proposition}\label{p:uniq_cont}
  Suppose that either $u$ is a weak solution to (\ref{prob1}) under
  the assumptions of Theorem \ref{asym1} or $u$ satisfies
  (\ref{eq:u1}--\ref{eq:u2}) and weakly solves (\ref{prob2}) under the
  assumptions of Theorem \ref{asym2}. If there exists $t_0\in(0,T)$
  such that
$$
u(x,t_0)=0\quad\text{for a.e. }x\in\R^N,
$$
then $u\equiv 0$ in $\R^N\times(0,T)$.
\end{Proposition}

There exists a large literature dealing with strong continuation
properties in the parabolic setting. \cite{LI1} (see too \cite{LI2})
studies parabolic operators with $L^{\frac{N+1}{2}}$ time-independent
coefficients obtaining a unique continuation property at a fixed time
$t_0$: the used technique relies on a representation formula for
solutions of parabolic equations in terms of eigenvalues of the
corresponding elliptic operator and cannot be applied to more general
equations with time-dependant coefficients.  \cite{SS} and \cite{SO}
use parabolic variants of the Carleman weighted inequalities to obtain
a unique continuation property at fixed time $t_0$ for parabolic
operators with time-dependant coefficients. In this direction, it is
worth mentioning the work of Chen \cite{CH} which contains not only a
unique continuation result but also some local asymptotic analysis of
solutions to parabolic inequalities with bounded coefficients; the
approach is based in recasting equations in terms of parabolic
self-similar variables. We also quote \cite{AV,E,EF,EKPV,EV,FE} for
unique continuation results for parabolic equations with
time-dependent potentials by Carleman inequalities and monotonicity
methods.

The present paper is organized as follows. In section
\ref{sec:parabolic-hardy-type}, we state some parabolic Hardy type
inequalities and weighted Sobolev embeddings related to  equations
(\ref{prob1}) and (\ref{prob2}). In section \ref{sec:spectr-analys-self}, we
completely describe the spectrum of the operator $L$ defined in
(\ref{eq:13}) and prove Proposition \ref{p:explicit_spectrum}. Section
\ref{sec:almgren} contains an Almgren parabolic monotonicity formula
which provides the unique continuation principle stated in Proposition
\ref{p:uniq_cont} and is used in section \ref{sec:blow-up-analysis},
together with a blow-up method, to prove Theorems \ref{asym1} and
 \ref{asym2}.

\medskip
\noindent
{\bf Notation. } We list below some notation used throughout the
paper.\par
\begin{itemize}
\item[-] ${\rm const}$ denotes some positive constant which may vary
  from formula to formula.
\item[-] $dS$ denotes the volume element on the unit
  $(N-1)$-dimensional sphere ${\mathbb S}^{N-1}$.
\item[-] $\omega_{N-1}$ denotes the volume of  ${\mathbb
  S}^{N-1}$, i.e. $\omega_{N-1}=\int_{{\mathbb S}^{N-1}}dS(\theta)$.
\item[-] For all $s\in\R$, $(s)_i$ denotes the Pochhammer's symbol
$(s)_i=\prod_{j=0}^{i-1}(s+j)$, $(s)_0=1$.

\end{itemize}

\section{Parabolic Hardy type
  inequalities and Weighted Sobolev embeddings}\label{sec:parabolic-hardy-type}
The following lemma provides a Hardy type inequality for parabolic
operators. We refer to \cite[Proposition 3.1]{poon} for a proof.

\begin{Lemma} \label{Hardytemp}
For every $t>0$ and $u\in {\mathcal H}_t$ there holds
$$
\int_{\R^N}\dfrac{u^2(x)}{|x|^2}\,G(x,t)\,dx\leq
\dfrac{1}{(N-2)t}\dint_{\R^N}u^2(x)G(x,t)\,dx+
\dfrac{4}{(N-2)^2}\dint_{\R^N}|\nabla
u(x)|^2\,G(x,t)\,dx.
$$
\end{Lemma}
\noindent
In the anisotropic version of the above inequality, a crucial role is
played by the first eigenvalue $\mu_1(a)$ of the angular
operator $-\Delta_{\mathbb S^{N-1}}-a(\theta)$ on the unit sphere
 $\mathbb S^{N-1}$ defined in (\ref{eq:67}).
\begin{Lemma}\label{Hardy_aniso}
  For every $a\in L^{\infty}\big({\mathbb S}^{N-1}\big)$, $t>0$, and
  $u\in {\mathcal H}_t$, there holds
\begin{multline*}
\int_{\R^N}\bigg(|\nabla
u(x)|^2-\frac{a(x/|x|)}{|x|^2}\,u^2(x)\bigg)
\,G(x,t)\,dx+\frac{N-2}{4t}\int_{\R^N}u^2(x)G(x,t)\,dx\\
\geq \bigg(\mu_1(a)+\frac{(N-2)^2}4\bigg)
\int_{\R^N}\dfrac{u^2(x)}{|x|^2}\,G(x,t)\,dx.
\end{multline*}
\end{Lemma}
\begin{pf}
Let $u\in C^\infty_{\rm
  c}(\R^N\setminus\{0\})$.  The gradient of $u$ can be
written in polar coordinates as
\begin{displaymath}
  \nabla u(x)=\big(\partial_ru(r,\theta)\big)\theta+
  \frac1r\nabla_{{\mathbb S}^{N-1}}u(r,\theta),
\quad r=|x|,\quad \theta=\frac{x}{|x|},
\end{displaymath}
hence
\begin{displaymath}
     |\nabla u(x)|^2
=\big|\partial_ru(r,\theta)\big|^2+
\frac1{r^2}\big|\nabla_{{\mathbb S}^{N-1}}u(r,\theta)
\big|^2
  \end{displaymath}
and
\begin{multline}\label{eq:1}
  \int_{\R^N}\bigg(|\nabla
  u(x)|^2-\frac{a(x/|x|)}{|x|^2}\,u^2(x)\bigg) \,G(x,t)\,dx\\=
  t^{-\frac N2}\int_{{\mathbb S}^{N-1}}\bigg(
  \int_0^{+\infty}r^{N-1}e^{-\frac{r^2}{4t}}
  |\partial_ru(r,\theta)|^2\,dr\bigg)\,dS(\theta)\\
  +t^{-\frac N2}
  \int_0^{+\infty}\frac{r^{N-1}e^{-\frac{r^2}{4t}}}{r^2}\bigg(\int_{{\mathbb
      S}^{N-1}}\left[|\nabla_{{\mathbb
        S}^{N-1}}u(r,\theta)|^2-a(\theta)
    |u(r,\theta)|^2\right]\,dS(\theta)\bigg)\,dr.
\end{multline}
For all $\theta\in{\mathbb S}^{N-1}$, let $\varphi_{\theta}\in
C^\infty_{\rm c}((0,+\infty))$ be defined by
$\varphi_{\theta}(r)=u(r,\theta)$, and $\widetilde\varphi_{\theta}\in
C^\infty_{\rm c}(\R^N\setminus\{0\})$ be the radially symmetric
function given by
$\widetilde\varphi_{\theta}(x)=\varphi_{\theta}(|x|)$.
From Lemma \ref{Hardytemp}, it follows that
\begin{align}\label{eq:2}
  t^{-\frac N2} \int_{{\mathbb S}^{N-1}}&\bigg(
  \int_0^{+\infty}r^{N-1}e^{-\frac{r^2}{4t}}
  |\partial_ru(r,\theta)|^2\,dr\bigg)\,dS(\theta)\\
  \notag& = t^{-\frac
    N2}\int_{{\mathbb S}^{N-1}}\bigg( \int_0^{+\infty}r^{N-1}
  e^{-\frac{r^2}{4t}}|\varphi_{\theta}'(r)|^2\,dr\bigg)\,dS(\theta)\\
  \notag& =\frac1{\omega_{N-1}} \int_{{\mathbb S}^{N-1}}\bigg(
  \int_{\R^N}|\nabla\widetilde
  \varphi_{\theta}(x)|^2G(x,t)\,dx\bigg)\,dS(\theta) \\
  &\notag\geq \frac1{\omega_{N-1}} \frac{(N-2)^2}{4} \int_{{\mathbb
      S}^{N-1}}\bigg(\int_{\R^N}\frac{|\widetilde
    \varphi_{\theta}(x)|^2}{|x|^2}G(x,t)\,dx\bigg)\,dS(\theta)\\
  &\notag \quad-\frac1{\omega_{N-1}}\frac{N-2}{4t}\int_{{\mathbb
      S}^{N-1}}\bigg(\int_{\R^N}|\widetilde
  \varphi_{\theta}(x)|^2G(x,t)\,dx\bigg)\,dS(\theta)
  \\
  \notag&=t^{-\frac N2}\frac{(N-2)^2}4 \int_{{\mathbb S}^{N-1}}\bigg(
  \int_0^{+\infty}\frac{r^{N-1}e^{-\frac{r^2}{4t}}}{r^2}|u(r,\theta)|^2\,dr\bigg)
  \,dS(\theta)\\
  \notag&-t^{-\frac N2}\frac{N-2}{4t}\int_{{\mathbb S}^{N-1}}\bigg(
  \int_0^{+\infty}r^{N-1}e^{-\frac{r^2}{4t}}|u(r,\theta)|^2\,dr\bigg)\,dS(\theta)\\
  \notag&
  =
\frac{(N-2)^2}4
\int_{\R^N}\dfrac{u^2(x)}{|x|^2}\,G(x,t)\,dx-
\frac{N-2}{4t}\int_{\R^N}u^2(x)G(x,t)\,dx,
\end{align}
where $\omega_{N-1}$ denotes the volume of the unit sphere ${\mathbb
  S}^{N-1}$, i.e. $\omega_{N-1}=\int_{{\mathbb S}^{N-1}}dS(\theta)$.
On the other hand, from the definition of $\mu_1(a)$ it follows
that
\begin{equation} \label{eq:3}
  \int_{{\mathbb
      S}^{N-1}}\!\!\left[|\nabla_{{\mathbb S}^{N-1}}u(r,\theta)|^2\!-
a(\theta)|u(r,\theta)|^2\right]dS(\theta)
  \geq \mu_1(a)\int_{{\mathbb S}^{N-1}}\!|u(r,\theta)|^2dS(\theta).
\end{equation}
From (\ref{eq:1}), (\ref{eq:2}), and (\ref{eq:3}), we deduce that
\begin{multline*}
\int_{\R^N}\bigg(|\nabla
u(x)|^2-\frac{a(x/|x|)}{|x|^2}\,u^2(x)\bigg)
\,G(x,t)\,dx+\frac{N-2}{4t}\int_{\R^N}u^2(x)G(x,t)\,dx\\
\geq \bigg(\mu_1(a)+\frac{(N-2)^2}4\bigg)
\int_{\R^N}\dfrac{u^2(x)}{|x|^2}\,G(x,t)\,dx,
\end{multline*}
for all $u\in C^\infty_{\rm
    c}(\R^N\setminus\{0\})$, thus yielding the required inequality
by  density of $C^\infty_{\rm c}(\R^N\setminus\{0\})$
in ${\mathcal H}_t$.~\end{pf}

The following corollary provides a norm in ${\mathcal H}_t$ equivalent
to $\|\cdot\|_{{\mathcal H}_t}$ and naturally related to the heat
operator with the Hardy potential of equation (\ref{prob}).
\begin{Corollary}\label{c:pos_def}
  Let $a\in L^{\infty}\big({\mathbb S}^{N-1}\big)$
satisfying (\ref{eq:posde}).
Then, for every $t>0$,
\begin{align*}
  &\inf_{u\in {\mathcal
      H}_t\setminus\{0\}}\frac{\int_{\R^N}\big(|\nabla
    u(x)|^2-\frac{a(x/|x|)}{|x|^2}\,u^2(x)\big)
    \,G(x,t)\,dx+\frac{N-2}{4t}\int_{\R^N}u^2(x)G(x,t)\,dx}
  {\int_{\R^N}|\nabla u(x)|^2
    \,G(x,t)\,dx+\frac{N-2}{4t}\int_{\R^N}u^2(x)G(x,t)\,dx}\\[5pt]
  &\quad=\inf_{v\in {\mathcal
      H}\setminus\{0\}}\frac{\int_{\R^N}\big(|\nabla
    v(x)|^2-\frac{a(x/|x|)}{|x|^2}\,v^2(x)\big)
    \,G(x,1)\,dx+\frac{N-2}{4}\int_{\R^N}v^2(x)G(x,1)\,dx}
  {\int_{\R^N}|\nabla v(x)|^2
    \,G(x,1)\,dx+\frac{N-2}{4}\int_{\R^N}v^2(x)G(x,1)\,dx}>0.
\end{align*}
\end{Corollary}
\begin{pf}
  The equality of the two infimum levels follows  by
  the change of variables $u(x)=v(x/\sqrt t)$.  To prove
  that they are strictly positive, we argue by contradiction and
  assume that for every $\e>0$ there exists $v_\e\in {\mathcal
    H}\setminus\{0\}$ such that
\begin{multline*}
\int_{\R^N}\bigg(|\nabla
v_\e(x)|^2-\frac{a(x/|x|)}{|x|^2}\,v_\e^2(x)\bigg)
\,G(x,1)\,dx+\frac{N-2}{4}\int_{\R^N}v_\e^2(x)G(x,1)\,dx\\
<\e \bigg(\int_{\R^N}|\nabla
v_\e(x)|^2\,G(x,1)\,dx+\frac{N-2}{4}\int_{\R^N}v_\e^2(x)G(x,1)\,dx\bigg),
\end{multline*}
which,  by Lemma \ref{Hardy_aniso}, implies that
\begin{multline*}
  \bigg(\mu_1\bigg(\frac{a}{1-\e}\bigg)+\frac{(N-2)^2}4\bigg)
  \int_{\R^N}\dfrac{v_\e^2(x)}{|x|^2}\,G(x,1)\,dx\\
  \leq \int_{\R^N}\bigg(|\nabla
  v_\e(x)|^2-\frac{a(x/|x|)}{(1-\e)|x|^2}\,v_\e^2(x)\bigg)
  \,G(x,1)\,dx+\frac{N-2}{4}\int_{\R^N}v_\e^2(x)G(x,1)\,dx<0
\end{multline*}
and consequently
$$
\mu_1\bigg(\frac{a}{1-\e}\bigg)+\frac{(N-2)^2}4<0.
$$
By continuity of the map $a\mapsto \mu_1(a)$ with respect to the
$L^{\infty}\big({\mathbb S}^{N-1}\big)$-norm, letting $\e\to0$ the
above inequality yields $\mu_1(a)+\frac{(N-2)^2}4\leq0$, giving rise to
a contradiction with (\ref{eq:posde}).
\end{pf}

The above results combined with the negligibility assumption
\eqref{eq:h} on $h$ allow estimating  the quadratic form associated to
the linearly perturbed equation (\ref{prob1}) for small times as follows.
\begin{Corollary}\label{c:pos_per}
  Let $a\in L^{\infty}\big({\mathbb S}^{N-1}\big)$ satisfy
  (\ref{eq:posde}) and $h\in L^\infty_{{\rm loc}}(\R^N\setminus
  \{0\}\times (0,T))$ satisfy (\ref{eq:h}). Then there exist $C_1',C_2>0$ and
  $\overline{T}_1>0$ such that for every $t\in(0,\overline{T}_1)$,
$s\in(0,T)$, and $u\in
  {\mathcal H}_t$ there holds
\begin{align*}
\int_{\R^N}\bigg(|\nabla
u(x)|^2&-\frac{a(x/|x|)}{|x|^2}\,u^2(x)-h(x,s)u^2(x)\bigg)
\,G(x,t)\,dx\\
&\geq C_1'\int_{\R^N}\frac{u^2(x)}{|x|^2}\,G(x,t)\,dx
-\frac{C_2}t\int_{\R^N}u^2(x)G(x,t)\,dx\\
\int_{\R^N}\bigg(|\nabla
u(x)|^2&-\frac{a(x/|x|)}{|x|^2}\,u^2(x)-h(x,s)u^2(x)\bigg)
\,G(x,t)\,dx
+\frac{N-2}{4t}\int_{\R^N}u^2(x)G(x,t)\,dx\\
&\geq C_1'\bigg(\int_{\R^N}|\nabla
u(x)|^2\,G(x,t)\,dx+\frac1t\int_{\R^N}u^2(x)G(x,t)\,dx\bigg).
\end{align*}
\end{Corollary}
\begin{pf}
From (\ref{eq:h}), we have that, for every $u\in
  {\mathcal H}_t$, there holds
\begin{align}\label{eq:11}
&  \left| \int_{\R^N}h(x,s)u^2(x)G(x,t)\,dx\right| \leq C_h\bigg(
  \int_{\R^N}u^2(x)G(x,t)\,dx+\int_{\R^N}|x|^{-2+\e}u^2(x)G(x,t)\,dx\bigg)\\
\notag  &\leq C_h\bigg(
  \int_{\R^N}u^2(x)G(x,t)\,dx+t^{\e/2}\int\limits_{|x|\leq\sqrt t
  }\frac{u^2(x)}{|x|^2}G(x,t)\,dx
+t^{-1+\e/2}\int\limits_{|x|\geq\sqrt t
  }u^2(x)G(x,t)\,dx
\bigg)\\
\notag&=\frac{C_h}{t}(t+t^{\e/2})\int_{\R^N}u^2(x)G(x,t)\,dx+
C_ht^{\e/2}\int_{\R^N}\frac{u^2(x)}{|x|^2}G(x,t)\,dx.
\end{align}
The stated inequalities follow from (\ref{eq:11}),
Lemma \ref{Hardytemp}, Corollary \ref{c:pos_def},
and assumption  (\ref{eq:posde}).
\end{pf}

In order to estimate the quadratic form associated to the nonlinearly
perturbed equation (\ref{prob2}), we derive a Sobolev type embedding in
spaces ${\mathcal H}_t$. To this purpose, we need the following
inequality, whose proof can be found in \cite[Lemma 3]{EFV}.

\begin{Lemma}\label{l:ineqx2}
For every $u\in {\mathcal H}$, $|x|u\in{\mathcal L}$ and
$$
\frac{1}{16}
\int_{\R^N}|x|^2u^2(x)G(x,1)\,dx\leq
\int_{\R^{N}}|\nabla u(x)|^2G(x,1)\,dx+\frac{N}{4}\int_{\R^N}u^2(x)G(x,1)\,dx.
$$
\end{Lemma}

The change of variables $u(x)=v(x/\sqrt t)$ in Lemma \ref{l:ineqx2},
yields the following inequality in ${\mathcal H}_t$.
\begin{Corollary}\label{cor:ineq}
For every $u\in {\mathcal H}_t$, there holds
$$
\frac{1}{16t^2}
\int_{\R^N}|x|^2u^2(x)G(x,t)\,dx\leq
\int_{\R^{N}}|\nabla u(x)|^2G(x,t)\,dx+\frac{N}{4t}\int_{\R^N}u^2(x)G(x,t)\,dx.
$$
\end{Corollary}

From Lemma \ref{l:ineqx2} and classical Sobolev embeddings, we can
easily deduce the following weighted Sobolev inequality (see also
\cite{ES}).

\begin{Lemma}\label{l:sob}
For all $u\in {\mathcal H}$ and $s\in[2,2^*]$, there holds
 $u\sqrt{G(\cdot,1)}\in L^{s}(\R^N)$. Moreover,
for every $s\in[2,2^*]$ there exists $C_s>0$ such that
$$
\bigg(
\int_{\R^N}|u(x)|^sG^{\frac{s}{2}}(x,1)\,dx\bigg)^{\!\!\frac{2}{s}}\leq C_s\bigg(
\int_{\R^{N}}\big(|\nabla u(x)|^2+u^2(x)\big)G(x,1)\,dx\bigg)
$$
for all $u\in {\mathcal H}$.
\end{Lemma}
\begin{pf}
  From Lemma \ref{l:ineqx2}, it follows that, if $u\in{\mathcal H}$,
  then $u\sqrt{G(\cdot,1)}\in H^1(\R^N)$; hence, by classical Sobolev
  embeddings, $u\sqrt{G(\cdot,1)}\in L^{s}(\R^N)$ for all
  $s\in[2,2^*]$. The stated inequality follows from classical Sobolev
  inequalities and Lemma \ref{l:ineqx2}.
\end{pf}

The change of variables $u(x)=v(x/\sqrt t)$ in Lemma \ref{l:sob},
yields the following inequality in ${\mathcal H}_t$.
\begin{Corollary}\label{cor:ineqSob}
  For every $t>0$, $u\in {\mathcal H}_t$, and  $2\leq s\leq 2^{*}$,
  there holds
$$
\Big(
\int_{\R^N}|u(x)|^{s}G^{\frac{s}{2}}(x,t)\,dx\Big)^{\!\!\frac{2}{s}}\leq C_s
t^{-\frac{N}{s}\left(\frac{s-2}{2}\right)}\|u\|^{2}_{{\mathcal H}_t}.
$$
\end{Corollary}

The above Sobolev estimate allows proving the nonlinear counterpart
of Corollary \ref{c:pos_per}.
\begin{Corollary}\label{c:pos_per_nonlin}
  Let $a\in L^{\infty}\big({\mathbb S}^{N-1}\big)$ satisfy
  (\ref{eq:posde}) and $\varphi\in C^1(\R^N\times(0,T)\times\R)$
  satisfy (\ref{eq:fi}) for some $1\leq p<2^{*}-1$. Then there
  exist $C_1''>0$ and a function $\overline{T}_2:(0,+\infty)\to\R$ such
  that, for every $R>0$, $t\in (0,\overline{T}_2(R))$, $s\in(0,T)$,
  and $u\in \{v\in {\mathcal H}_t\cap L^{p+1}(\R^N):
  \|v\|_{L^{p+1}(\R^N)}\leq R\}$, there holds
\begin{align*}
\int_{\R^N}\bigg(|\nabla
u(x)|^2&-\frac{a(x/|x|)}{|x|^2}\,u^2(x)-\varphi(x,s,u(x))u(x)\bigg)
\,G(x,t)\,dx
+\frac{N-2}{4t}\int_{\R^N}u^2(x)G(x,t)\,dx\\
&\geq C_1''\bigg(\int_{\R^N}|\nabla
u(x)|^2\,G(x,t)\,dx+\frac1t\int_{\R^N}u^2(x)G(x,t)\,dx\bigg).
\end{align*}
\end{Corollary}
\begin{pf}
  From (\ref{eq:fi}),  H\"{o}lder's inequality, and
  Corollary \ref{cor:ineqSob}, we have that, for all $u\in {\mathcal
    H}_t\cap L^{p+1}(\R^N)$, there holds
\begin{align}\label{eq:varphi}
   \bigg| \int_{\R^N}&\varphi(x,s,u(x))u(x)G(x,t)\,dx\bigg| \\
&\notag\leq
  C_\varphi\bigg(
  \int_{\R^N}\!\!\!u^{2}(x)G(x,t)\,dx+\int_{\R^N}\!\!\!
  u^{2}(x)|u(x)|^{p-1}G(x,t)\,dx\bigg)\\
  \notag &\leq C_{\varphi}\bigg( \int_{\R^N}\!\!\!u^2(x)G(x,t)\,dx+
\bigg(\int_{\R^N}|u(x)|^{p+1}G^{\frac{p+1}{2}}(x,t)\,dx
  \bigg)^{\!\!\frac{2}{p+1}}\|u\|^{p-1}_{L^{p+1}(\R^N)}\bigg)\\
  \notag &\leq C_{\varphi}\bigg(
C_{p+1}  t^{\frac{(N+2)-p(N-2)}{2(p+1)}}\|u\|^{p-1}_{L^{p+1}(\R^N)}
  \int_{\R^N}|\nabla u(x)|^2\,G(x,t)\,dx\\
  \notag&\hskip2cm +
  \Big(t+C_{p+1}t^{\frac{(N+2)-p(N-2)}{2(p+1)}}\|u\|^{p-1}_{L^{p+1}(\R^N)}\Big)
  \frac1t\int_{\R^N}u^2(x)G(x,t)\,dx\bigg)
\end{align}
with $C_{p+1}$ as in Corollary \ref{cor:ineqSob}.
The stated inequality follows from Corollary \ref{c:pos_def} and
(\ref{eq:varphi}) by choosing $t$ sufficiently small depending on
$\|u\|_{L^{p+1}(\R^N)}$.
\end{pf}

\section{Spectrum of Ornstein-Uhlenbeck type operators with inverse
  square potentials}\label{sec:spectr-analys-self}

In this section we describe the spectral properties of the operator
$L$ defined in (\ref{eq:13}),
extending to anisotropic singular potentials the analysis carried out
in \cite{vazquez_zuazua} for $a\equiv\lambda$ constant.
Following \cite{ES}, we first prove the following compact embedding.

\goodbreak
\begin{Lemma}\label{l:compact}
The space ${\mathcal H}$ is compactly embedded in
${\mathcal L}$.
\end{Lemma}
\begin{pf}
Let us  assume that $u_{k}\weakly
u$ weakly in ${\mathcal H}$. From Rellich's theorem
$u_{k}\rightarrow u$ in  $L^{2}_{\rm loc}(\R^N)$. For every $R>0$ and $k\in\N$,
we have
\begin{equation}\label{eq:6}
\dint_{\R^N}|u_{k}-u|^2G(x,1)\,dx=
A_k(R)+B_k(R)
\end{equation}
where
\begin{equation}\label{eq:7}
A_k(R)=\int_{\{|x|\leq R\}}|u_{k}(x)-u(x)|^2e^{-|x|^2/{4}}\,dx
\to 0\quad\text{as }k\to+\infty,\quad \text{for every }R>0
\end{equation}
and
$$
B_k(R)=
\int_{\{|x|>R\}}|u_{k}(x)-u(x)|^2G(x,1)\,dx.
$$
From Lemma \ref{l:ineqx2} and boundedness of $u_k$ in ${\mathcal H}$,
we deduce that
\begin{align}\label{eq:8}
&B_{k}(R)\leq
R^{-2}\dint_{\{|x|>R\}}|x|^2|u_{k}(x)-u(x)|^2G(x,1)\,dx\\
\notag&\leq
\frac{1}{R^2}\bigg(16\int_{\R^N}|\nabla
(u_{k}-u)(x)|^2G(x,1)\,dx+4N\int_{\R^N}|u_{k}(x)-u(x)|^2G(x,1)\,dx\bigg)
\leq \frac{\rm const}{R^2}.
\end{align}
Combining (\ref{eq:6}), (\ref{eq:7}), and (\ref{eq:8}), we obtain that
 $u_{k}\rightarrow u$ strongly in ${\mathcal L}$.
\end{pf}

From classical spectral theory we deduce the following abstract
description of the spectrum of~$L$.
\begin{Lemma}\label{l:hilbasis}
  Let $a\in L^{\infty}\big({\mathbb S}^{N-1}\big)$ such that
  (\ref{eq:posde}) holds. Then the spectrum of the operator $L$
  defined in (\ref{eq:13}) consists of a diverging sequence of real
  eigenvalues with finite multiplicity. Moreover, there exists an
  orthonormal basis of ${\mathcal L}$ whose elements belong to
  ${\mathcal H}$ and are eigenfunctions of $L$.
\end{Lemma}
\begin{pf}
  By Corollary \ref{c:pos_def} and the Lax-Milgram Theorem, the bounded
  linear self-adjoint operator
$$
T:{\mathcal L}\to{\mathcal L},\quad
T=\bigg(L+\frac{N-2}{4}\,{\rm Id}\bigg)^{-1}
$$
is well defined. Moreover, by Lemma \ref{l:compact}, $T$ is compact.
The result then follows from the Spectral Theorem.
\end{pf}

Let us now compute explicitly the eigenvalues of $L$ with the
corresponding multiplicities and eigenfunctions by proving Proposition
\ref{p:explicit_spectrum}.

\begin{pfn}{Proposition \ref{p:explicit_spectrum}}
  Assume that $\gamma$ is an eigenvalue of $L$ and $g\in{\mathcal
    H}\setminus\{0\}$ is a corresponding eigenfunction, so that
\begin{equation}\label{eq:50}
-\Delta g(x)+ \frac{\nabla g(x)\cdot
  x}{2}-\frac{a(x/|x|)}{|x|^2}\,g(x)=\gamma\, g(x)
\end{equation}
in a weak ${\mathcal H}$-sense.  From classical regularity theory for
elliptic equations, $g\in C^{1,\alpha}_{\rm loc}(\R^N\setminus\{0\})$.
 Hence $g$ can be
expanded as
\begin{equation*}
g(x)=g(r\theta)=\sum_{k=1}^\infty\phi_k(r)\psi_k(\theta)
\quad \text{in }L^2({\mathbb S}^{N-1}),
\end{equation*}
where $r=|x|\in(0,+\infty)$, $\theta=x/|x|\in{{\mathbb S}^{N-1}}$, and
\begin{equation*}
  \phi_k(r)=\int_{{\mathbb S}^{N-1}}g(r\theta)
  \psi_k(\theta)\,dS(\theta).
\end{equation*}
Equations (\ref{eq:2rad}) and (\ref{eq:50}) imply that, for every $k$,
\begin{equation}\label{eq:51}
  \phi''_{k}+\left(\dfrac{N-1}{r}-\dfrac{r}{2}\right)
  \phi'_{k}+\left(\gamma-\dfrac{\mu_k}{r^2}\right)\phi_{k}=0
\quad\text{in  }(0,+\infty).
\end{equation}
Since  $g\in {\mathcal H}$, we have that
\begin{equation}\label{condition1}
\infty>\int_{\R^N}g^2(x)G(x,1)\,dx=\int_{0}^{\iy}
\!\bigg(\int_{{\mathbb S}^{N-1}}g^{2}(r\theta)\,dS(\theta)\bigg)
r^{N-1}e^{-\frac{r^2}{4}}\,dr\geq \int_{0}^{\iy}r^{N-1}e^{-\frac{r^2}{4}}\phi_{k}^{2}(r)\,dr
\end{equation}
and, by  the Hardy type inequality of Lemma \ref{Hardytemp},
\begin{equation}\label{condition2}
\infty>\int_{\R^N}\dfrac{g^2(x)}{|x|^2}
G(x,1)\,dx\geq\int_{0}^{\iy}r^{N-3}e^{-\frac{r^2}{4}}\phi_{k}^2(r)\,dr.
\end{equation}
For all $k=1,2,\dots$ and $t>0$, we define
$w_{k}(t)=(4t)^{\frac{\alpha_k}{2}}\phi_k(\sqrt{4t})$, with
$\alpha_{k}=\frac{N-2}{2}-\sqrt{\big(\frac{N-2}2\big)^{\!2}+\mu_{k}(a)}$.
From (\ref{eq:51}), $w_k$ satisfies
\begin{equation*}
  t w_{k}''(t)+\left(\frac{N}{2}-\alpha_k-t\right)w'_{k}(t)+
  \left(\frac{\alpha_k}{2}+\gamma\right)w_{k}(t)=0\quad\text{in }(0,+\infty).
\end{equation*}
Therefore, $w_{k}$ is a solution of the well known Kummer Confluent
Hypergeometric Equation (see \cite{Abramowitz_Stegun} and
\cite{macdonald}). Then there exist $A_k,B_k\in\R$ such that
$$
w_k(t)=A_k M\Big(-\frac{\alpha_k}{2}-\gamma,\frac N2-\alpha_k,t\Big)
+B_k U\Big(-\frac{\alpha_k}{2}-\gamma,\frac N2-\alpha_k,t\Big),
\quad t\in (0,+\infty).
$$
Here $M(c,b,t)$ and, respectively, $U(c,b,t)$ denote the Kummer
function (or confluent hypergeometric function) and, respectively, the
Tricomi function (or confluent hypergeometric function of the second
kind); $M(c,b,t)$ and $U(c,b,t)$ are two linearly independent
solutions to the Kummer Confluent Hypergeometric Equation
$$
tw''(t)+(b-t)w'(t)-ct=0,\quad t\in (0,+\infty).
$$
Since $\big(\frac N2-\alpha_k\big)>1$, from the well-known asymptotics
of $U$ at $0$ (see e.g. \cite{Abramowitz_Stegun}), we have that
$$
U\Big(-\frac{\alpha_k}{2}-\gamma,\frac N2-\alpha_k,t\Big)
\sim \text{\rm const}\,t^{1-\frac{N}{2}+\alpha_k}
\quad\text{as }t\to 0^+,
$$
for some $\text{\rm const}\neq 0$ depending only on
$N,\gamma$, and $\alpha_k$. On the other hand, $M$ is the sum of the series
$$
M(c,b,t)=\sum_{n=0}^\infty
\frac{(c)_n}{(b)_n}\,\frac{t^n}{n!}.
$$
We notice that $M$ has a finite limit at $0^+$, while its behavior at
$\infty$ is singular and depends on the value
$-c=\frac{\alpha_k}{2}+\gamma$.  If $\frac{\alpha_k}{2}+\gamma=m\in
\N=\{0,1,2,\cdots\}$, then $ M\big(-\frac{\alpha_k}{2}-\gamma,\frac
N2-\alpha_k,t\big)$ is a polynomial of degree $m$ in $t$, which we
will denote as $P_{k,m}$, i.e.,
$$
P_{k,m}(t)=M\Big(-m,{\textstyle{\frac N2}}-\alpha_k,t\Big)=
\sum_{n=0}^m
\frac{(-m)_n}{\big(\frac{N}2-\alpha_k\big)_n}\,\frac{t^n}{n!}.
$$
If $\big(\frac{\alpha_k}{2}+\gamma\big)\not\in \N$, then from the well-known
asymptotics of $M$ at $\infty$
(see e.g. \cite{Abramowitz_Stegun})
we have that
$$
M\Big(-\frac{\alpha_k}{2}-\gamma,\frac N2-\alpha_k,t\Big)
\sim \text{\rm const}\,e^tt^{-\frac{N}{2}+\frac{\alpha_k}2-\gamma}
\quad\text{as }t\to +\infty,
$$
for some $\text{\rm const}\neq 0$ depending only on
$N,\gamma$, and $\alpha_k$.

Now, let us fix $k\in\N$, $k\geq 1$. From the above description, we have that
$$
w_k(t)\sim {\rm const\,}B_k
t^{1-\frac{N}{2}+\alpha_k}
\quad\text{as }t\to 0^+,
$$
for some $\text{\rm const}\neq 0$,
and hence
$$
\phi_k(r)=r^{-\alpha_k}w_k\Big(\frac{r^2}{4}\Big)\sim
{\rm const\,}B_k
r^{2-N+\alpha_k}
\quad\text{as }r\to 0^+,
$$
for some $\text{\rm const}\neq 0$. Therefore, condition (\ref{condition2})
can be satisfied only for $B_k=0$.
If $\frac{\alpha_k}{2}+\gamma\not\in \N$,
then
$$
w_k(t)\sim {\rm const\,}A_ke^t
t^{-\frac{N}{2}+\frac{\alpha_k}2-\gamma}
\quad\text{as }t\to +\infty,
$$
for some $\text{\rm const}\neq 0$,
and hence
$$
\phi_k(r)=r^{-\alpha_k}w_k\Big(\frac{r^2}{4}\Big)\sim
{\rm const\,}A_k
r^{-N-2\gamma}e^{r^2/4}
\quad\text{as }r\to +\infty,
$$
for some $\text{\rm const}\neq 0$. Therefore, condition (\ref{condition1})
can be satisfied only for $A_k=0$.
 If $\frac{\alpha_k}{2}+\gamma=m\in \N$, then
$r^{-\alpha_k}P_{k,m}\big(\frac{r^2}{4}\big)$ solves (\ref{eq:51}); moreover
the function
$$
|x|^{-\alpha_k}P_{k,m}\Big(\frac{|x|^2}{4}\Big)
\psi_k\Big(\frac{x}{|x|}\Big)
$$
belongs to ${\mathcal H}$, thus  providing an eigenfunction of $L$.

We can conclude from the above discussion that if
$\frac{\alpha_k}{2}+\gamma\not\in \N$ for all $k\in\N$, $k\geq 1$,
then $\gamma$ is not an eigenvalue of $L$. On the other hand,
if there exist $k_0,m_0\in\N$, $k_0\geq1$, such that
$$
\gamma=\gamma_{m_0,k_0}=m_0-\frac{\alpha_{k_0}}{2}
$$
then $\gamma$ is an eigenvalue of $L$
with multiplicity
\begin{equation}\label{eq:cardinal}
m(\gamma)=m(\gamma_{m_0,k_0})=\#\bigg\{j\in\N,j\geq 1:
\gamma_{m_0,k_0}+\frac{\alpha_j}2\in\N\bigg\}<+\infty
\end{equation}
and a basis of the corresponding eigenspace is
$$
\left\{
  |x|^{-\alpha_j}P_{j,\gamma_{m_0,k_0}+{\alpha_j}/2}\bigg(\frac{|x|^2}{4}\bigg)
  \psi_j\Big(\frac{x}{|x|}\Big): j\in\N,j\geq
  1,\gamma_{m_0,k_0}+\frac{\alpha_j}2\in\N \right\}.
$$
The proof is thereby complete.
\end{pfn}

\begin{remark}\label{rem:chen}
  If $a(\theta)\equiv0$, then $\mu_k(0)=k(N+k-2)$, so that
  $\alpha_k=\frac{(N-2)}{2}-\sqrt{\big(\frac{N-2}{2}+k\big)^2}=-k$,
  and $\gamma_{m,k}=\frac{k}{2}+m$. Hence, in this case we recover the
  well known fact (see e.g. \cite{CH} and \cite{ES}) that the
  eigenvalues of the Ornstein-Uhlenbeck operator
  $-\Delta+\frac{x}2\cdot\nabla$ are the positive half-integer
  numbers.
\end{remark}

\begin{remark}\label{rem:ortho}
  Due to orthogonality of eigenfunctions $\{\psi_k\}_k$ in
  $L^2({\mathbb S}^{N-1})$, it is easy to verify that
$$
\text{if }(m_1,k_1)\neq(m_2,k_2)\quad\text{then}\quad
V_{m_1,k_1}\text{ and } V_{m_2,k_2}\text{ are orthogonal in }{\mathcal
  L}.
$$
By Lemma \ref{l:hilbasis}, it follows that
$$
\left\{
\widetilde V_{n,j}=
\frac{V_{n,j}}{\|V_{n,j}\|_{\mathcal L}}: j,n\in\N,j\geq
  1\right\}
$$
is an orthonormal basis of ${\mathcal L}$.
\end{remark}

\section{The parabolic Almgren monotonicity formula }\label{sec:almgren}

Throughout this section, we will assume that $a\in
L^{\infty}\big({\mathbb S}^{N-1}\big)$ satisfies (\ref{eq:posde}) and
\emph{either}
\begin{align}\label{eq:82}\tag{${\bf I}$}
\text{$u$ is a weak solution to (\ref{prob1})  with  $h$ satisfying
  (\ref{eq:der}) and (\ref{eq:h})}
\end{align}
 \emph{or}
\begin{align}\label{eq:83}\tag{${\bf II}$}
  \text{$u$ satisfies (\ref{eq:u1}--\ref{eq:u2}) and weakly solves
    (\ref{prob2}) for some $\varphi\in C^1(\R^N\times(0,T)\times\R)$
    satisfying (\ref{eq:fi})}.
\end{align}
We denote as
$$
f(x,t, s)=
\begin{cases}
h(x,t)s,&\text{ in case {\bf (I)}},\\
\varphi(x,t, s),&\text{ in case {\bf (II)}},
\end{cases}
$$
so that, in both cases, $u$ is a weak solution to (\ref{prob}) in
$\R^N\times(0,T)$ in the sense of Definition \ref{def:solution}.  Let
\begin{equation}\label{eq:5}
\overline T=
\begin{cases}
\overline{T}_1,&\text{ in case {\bf (I)}},\\
\overline{T}_2(R_0),&\text{ in case {\bf (II)}},
\end{cases}
\quad\text{and}\quad
C_1=
\begin{cases}
C_1',&\text{ in case {\bf (I)}},\\
C_1'',&\text{ in case {\bf (II)}},
\end{cases}
\end{equation}
being $C_1',\overline{T}_1$ as in Corollary \ref{c:pos_per} and
$C_1'',\overline{T}_2(R_0)$  as in Corollary \ref{c:pos_per_nonlin} with
$$
R_0=\sup_{t\in(0,T)}\|u(\cdot,t)\|_{L^{p+1}(\R^N)}
$$
(notice that $R_0$ is finite by assumption (\ref{eq:u1})).  We denote
\begin{equation*}
  \alpha=\frac{T}{2\big(\big\lfloor{T}/{\overline T}\big\rfloor+1\big)},
\end{equation*}
where $\lfloor \cdot\rfloor$ denotes the floor function, i.e. $\lfloor
x\rfloor:=\max\{j\in\Z:\ j\leq x\}$.
Then
$$
(0,T)=\bigcup_{i=1}^k(a_i,b_i)
$$
where
$$
k=2\big(\big\lfloor{T}/{\overline T}\big\rfloor+1\big)-1,
\quad
a_i=(i-1)\alpha,\quad\text{and}\quad
b_i=(i+1)\alpha.
$$
We notice that $0<2\alpha<{\overline T}$ and $(a_i,b_i)\cap(a_{i+1},b_{i+1})=
(i\alpha,(i+1)\alpha)\not =\emptyset$.
For every $i$, $1\leq i\leq k$, we define
\begin{equation}\label{eq:76}
u_{i}(x,t)=u(x, t+a_i),\quad
x\in\R^N,\ t\in(0,2\alpha).
\end{equation}

\begin{Lemma}\label{l:u_i}
  For every $i=1,\dots,k$, the function $u_i$ defined in (\ref{eq:76})
  is a weak solution to
\begin{equation}\label{prob_i}
(u_i)_t+\Delta u_i+\dfrac{a(x/|x|)}{|x|^2}\,u_i+f(x,t+a_i, u_i(x,t))=0
\end{equation}
  in $\R^N\times (0,2\alpha)$ in
  the sense of Definition \ref{def:solution}.  Furthermore, the
  function $v_i(x,t):=u_i(\sqrt{t}x,t)$ is a weak solution to
\begin{equation}\label{eq:eqforv_i}
  (v_i)_t+\frac1t\bigg(\Delta v_i-\frac x2\cdot \nabla v_i+
\dfrac{a(x/|x|)}{|x|^2}\,v_i+tf\big(\sqrt tx,t+a_i,v_i(x,t)\big)\bigg)=0
\end{equation}
 in $\R^N\times(0,2\alpha)$ in the sense of Remark
  \ref{rem:v2}.
\end{Lemma}
\begin{pf} If $i=1$, then $a_1=0$, $u_1(x,t)=u (x,t)$ in $\R^N\times
  (0,2\alpha)$, and we immediately conclude.  For every $1<i\leq k$,
  $a_i\neq 0$, and, being $G(x,t)$ as in (\ref {eq:heatker}), the
  following properties hold for all $t\in(a_i,b_i)$:
\begin{align*}
{\rm (i)}\quad &G\big(x,{\textstyle{\frac{t(t-a_i)}{a_i}}}\big)G(x,t)=
\big({\textstyle{\frac{t^2}{a_i}}}\big)^{-N/2}G(x,t-a_i);\\
{\rm (ii)}\quad & \text{if $\phi\in {\mathcal H}_{t-a_i}$, then
  $\phi\,G\big(\cdot,{\textstyle{\frac{t(t-a_i)}{a_i}}}\big)
  \in {\mathcal H}_{t}$};\\
{\rm (iii)}\quad & \text{if $\psi\in ({\mathcal H}_t)^\star$, then
  $\psi\in
  ({\mathcal H}_{t-a_{i}})^\star$ and}\\
&{\phantom{\big\langle}}_{{\mathcal H}_{t-a_i}^\star}\big\langle
\psi,\phi \big\rangle_{{\mathcal
    H}_{t-a_i}}=\bigg(\dfrac{t^2}{a_i}\bigg)^{\!\!\frac{N}{2}}
\!\!\!\!\!{\phantom{\big\langle}}_{{\mathcal H}_{t}^\star}\Big\langle \psi,
\phi\,G\big(\cdot,{\textstyle{\frac{t(t-a_i)}{a_i}}}\big)
\Big\rangle_{{\mathcal H}_{t}},\quad \text{for all $\phi\in {\mathcal
    H}_{t-a_i}$}.
\end{align*}
Let $1<i\leq k$ and $\phi\in {\mathcal H}_{t-a_i}$. Due to (ii),
$\phi\,G\big(\cdot,\frac{t(t-a_i)}{a_i}\big)\in {\mathcal H}_{t}$ and
then, since $u$ is a solution to (\ref{prob}) in the sense of of
Definition \ref{def:solution}, for a.e. $t\in(a_i,b_i)$ we have
\begin{multline}
  {\phantom{\bigg\langle}}_{{\mathcal
      H}_t^\star}\bigg\langle u_t+\frac{\nabla u\cdot
    x}{2t},\phi\,G\big(x,{\textstyle{\frac{t(t-a_i)}{a_i}}}\big)
\bigg\rangle_{{\mathcal H}_t}\\= \int_{\R^N}\nabla u(x,t)\cdot \nabla
  \phi(x)\, G\big(x,{\textstyle{\frac{t(t-a_i)}{a_i}}}\big)G(x,t)\,dx
-
  \int_{\R^N}\phi(x) \dfrac{a_{i}x\cdot \nabla u(x,t)}{2(t-a_{i})t}
G\big(x,{\textstyle{\frac{t(t-a_i)}{a_i}}}\big)G(x,t)\,dx\\-
\int_{\R^N}\frac{a(x/|x|)}{|x|^2}\,u(x,t)\phi(x)
G\big(x,{\textstyle{\frac{t(t-a_i)}{a_i}}}\big)G(x,t)\,dx
 -
  \int_{\R^N}f(x,t,u(x,t))\phi(x)
G\big(x,{\textstyle{\frac{t(t-a_i)}{a_i}}}\big)G(x,t)\,dx.
\end{multline}
Therefore, thanks to (i) and (iii), we obtain
\begin{align*} {\phantom{\bigg\langle}}_{{\mathcal
      H}_{t-a_{i}}^\star}\bigg\langle u_t+\frac{\nabla u(x,t)\cdot
    x}{2(t-a_{i})} , \phi \bigg\rangle_{{\mathcal H}_{t-a_{i}}}=&
  \int_{\R^N}\bigg( \nabla u(x,t)\cdot \nabla \phi(x)
  -\frac{a(x/|x|)}{|x|^2}\,u(x,t)\phi(x)\bigg)\,G(x,t-a_i)\,dx\\
  &- \int_{\R^N}f(x,t,u(x,t))\phi(x) G(x,t-a_i)\,dx.
\end{align*}
By the change of variables $s=t-a_{i}$, we conclude that
$u_{i}(x,t)=u(x, t+a_i)$ is a weak solution to (\ref{prob_i}) in
$\R^N\times (0,2\alpha)$ in the sense of Definition
\ref{def:solution}. By a further change of variables, we easily obtain that
 $v_i(x,t):=u_i(\sqrt{t}x,t)$ is a weak
solution to (\ref{eq:eqforv_i}) in $\R^N\times(0,2\alpha)$ in the sense
of Remark~\ref{rem:v2}.~\end{pf}

\noindent
For every $i=1,\dots,k$, we  define
\begin{equation}\label{eq:Hi(t)}
  H_i(t)=\int_{\R^N}u_{i}^2(x,t)\, G(x,t)\,dx,
  \quad\text{for every }t\in (0,2\alpha),
\end{equation}
and
\begin{equation}\label{eq:Di(t)}
  D_i(t)=\!\!\int_{\R^N}\!\!\bigg(|\nabla u_i(x,t)|^2-
  \dfrac{a\big(\frac{x}{|x|}\big)}{|x|^2}u_i^2(x,t)-
f(x,t+a_i,u_i(x,t))u_{i}(x,t)\bigg)G(x,t)\,dx
\end{equation}
for a.e. $t\in (0,2\alpha)$.
\begin{Lemma}\label{l:Hprime}
For every $1\leq i\leq k$,
 $H_i\in W^{1,1}_{\rm loc}(0,2\alpha)$ and
\begin{equation}\label{eq:10i}
  H'_i(t)=2\!\!\!\!{\phantom{\bigg\langle}}_{{\mathcal
      H}_t^\star}\bigg\langle
  (u_i)_t+\frac{\nabla u_i\cdot x}{2t},u_i(\cdot,t)
  \bigg\rangle_{{\mathcal H}_t}=2D_i(t)\quad\text{for a.e. }t\in(0,2\alpha).
\end{equation}
\end{Lemma}
\begin{pf}
  It follows from Lemma \ref{l:u_i}  and Remark
  \ref{rem:uv}.
\end{pf}

\noindent
\begin{Lemma}\label{l:Hcreas}
  If  $C_1$ is as in (\ref{eq:5}), then, for every $i=1,\dots,k$, the
  function
$$
t\mapsto
 t^{-2C_1+\frac{N-2}{2}}H_i(t)$$
 is
nondecreasing in $(0,2\alpha)$.
\end{Lemma}
\begin{pf}
  From Lemma \ref{l:Hprime} and Corollaries \ref{c:pos_per} and
  \ref{c:pos_per_nonlin}, taking into account that $2\alpha<{\overline
    T}$, we have that, for all $t\in(0,2\alpha)$,
$$
H'_{i}(t)\geq \frac1t\bigg(2C_1-\frac{N-2}2\bigg)H_{i}(t),
$$
which implies
$$
\frac{d}{dt}\bigg(t^{-2C_1+\frac{N-2}{2}}H_{i}(t)\bigg)\geq 0.
$$
Hence the function $t\mapsto t^{-2C_1+\frac{N-2}{2}}H_i(t)$ is
nondecreasing in $(0,2\alpha)$.
\end{pf}

\noindent
\begin{Lemma}\label{l:Hpos}
If  $1\leq i\leq k$
and  $H_i(\bar t)=0$
  for some $\bar t\in(0,2\alpha)$, then $H_{i}(t)=0$ for all $t\in
  (0,\bar t\,]$.
\end{Lemma}
\begin{pf}
  From Lemma \ref{l:Hcreas}, the function $t\mapsto
  t^{-2C_1+\frac{N-2}{2}}H_{i}(t)$ is nondecreasing in $(0,2\alpha)$,
  nonnegative, and vanishing at $\bar t$. It follows that $H_i(t)=0$
  for all $t\in (0,\bar t]$.
\end{pf}

\noindent
The regularity of $D_i$ in $(0,2\alpha)$ is analyzed in the following
lemma.
\begin{Lemma}\label{l:Dprime}
  If $1\leq i\leq k$ and $T_i\in(0,2\alpha)$ is such that
  $u_i(\cdot,T_i)\in{\mathcal H}_{T_i}$, then
\begin{itemize}
\item[(i)]
$\dint_\tau^{T_i}\dint_{\R^N}\bigg(\left|(u_i)_t(x,t)+\frac{\nabla
      u_i(x,t)\cdot x}{2t}\right|^2G(x,t)\,dx\bigg)\,dt<+\infty
\quad\text{for all
  }\tau\in (0,T_i)$;\\[5pt]
\item[(ii)]
the function
\begin{equation*}
t\mapsto t D_i(t)
\end{equation*}
belongs to $W^{1,1}_{\rm loc}(0,T_i)$ and its weak derivative is,
for a.e. $t\in(0,T_i)$,  as follows:
\end{itemize}

\smallskip\noindent
{in case \bf (I)}
\begin{multline*}
  \frac{d}{dt}\,
  \big(tD_i(t)\big)=2t\int_{\R^N}\left|(u_i)_t(x,t)+\frac{\nabla
      u_i(x,t)\cdot x}{2t}\right|^2G(x,t)\,dx\\
  + \int_{\R^N}h(x,t+a_i)\left(\frac{N-2}{2}u_i^2(x,t)+(\nabla
    u_i(x,t)\cdot x)u_i(x,t)-\frac{|x|^2}{4t}u_i^2(x,t)\right)
  \,G(x,t)\,dx\\
  -t\int_{\R^N}h_{t}(x,t+a_i)u_i^2(x,t)G(x,t)\,dx;
\end{multline*}
{in  case \bf (II)}
\begin{align*}
  \frac{d}{dt}&\, \big(tD_i(t)\big)= 2t
  \int_{\R^N}\left|(u_i)_t(x,t)+\frac{\nabla
      u_i(x,t)\cdot x}{2t}\right|^2G(x,t)\,dx\\
  & +t\int_{\R^N} \bigg(
  \varphi(x,t+a_i,u_i(x,t))-\frac{\partial\varphi}{\partial u_i}
  (x,t+a_i,u_i(x,t))u_i(x,t)\bigg)(u_i)_t(x,t)G(x,t)\,dx\\
  &+\int_{\R^N} \bigg(\frac{N-2}2\varphi(x,t+a_i,u_i(x,t))u_i(x,t)-t
  \frac{\partial\varphi}{\partial t}(x,t+a_i,u_i(x,t))
  u_i(x,t)\\
&\hskip2cm-N\Phi(x,t+a_i,u_i(x,t))-\nabla_x\Phi(x,t+a_i,u_i(x,t))\cdot
  x\bigg)G(x,t)\,dx\\
  &+\int_{\R^N}\frac{|x|^2}{4t} \bigg(2\Phi(x,t+a_i,u_i(x,t))-
  \varphi(x,t+a_i,u_i(x,t))u_i(x,t)\bigg)G(x,t)\,dx
\end{align*}
where
$$
\Phi(x,t,s)=\int_0^s\varphi(x,t,\xi)\,d\xi.
$$
\end{Lemma}
\begin{pf} Let us first consider case {\bf (I)}, i.e.  $f (x,t,
  u)=h(x,t)u$, with $h(x,t)$ under conditions
  (\ref{eq:der}--\ref{eq:h}).  We test equation (\ref{eq:eqforv_i})
  with $(v_i)_t$; we notice that this is not an admissible test
  function for equation (\ref{eq:eqforv_i}) since a priori $(v_i)_t$
  does not take values in ${\mathcal H}$. However the formal testing
  procedure can be made rigorous by a suitable approximation.  Such a
  test combined with Corollary \ref{c:pos_per} yields, for all
  $t\in(0,T_i)$,
\begin{align*}
&  \int_t^{T_i}s\bigg(\int_{\R^N}(v_i)_t^2(x,s)G(x,1)\,dx\bigg)\,ds\leq{\rm
    const\,}\bigg( \|u_i(\sqrt {T_i}\,\cdot,T_i)\|^2_{\mathcal H}+
  \int_{\R^N}v_i^2(x,t)G(x,1)\,dx\\
  & +\int_t^{T_i}\!\!\bigg(\int_{\R^N}\!\!h(\sqrt s x, s+a_i)\bigg(
  \frac{|x|^2}8v_i^2(x,s)- \frac{\nabla v_i(x,s)\cdot x }{2}v_i(x,s)
  -\frac{N-2}{4}v_i^2(x,s)\bigg)G(x,1)\,dx\bigg)\,ds\\
  &+\frac{1}{2}\int_t^{T_i}\!\!s\bigg(\int_{\R^N}h_{s}(\sqrt s x,
  s+a_i)v_i^2(x,s)G(x,1)\,dx\bigg)\,ds\bigg).
\end{align*}
Since, in view of (\ref{eq:der}--\ref{eq:h}) and Lemmas
\ref{Hardytemp} and \ref{l:ineqx2}, the integrals in the last two
lines of the previous formula are finite for every $t\in(0,T_i)$, we
conclude that
$$
(v_i)_t\in L^2(\tau,T_i;{\mathcal L})\quad\text{for all
  }\tau\in (0,T_i).
$$
Testing (\ref{eq:eqforv_i}) with $(v_i)_t$ also yields
\begin{align*}
  &\int_t^{T_i}s\bigg(\int_{\R^N}(v_i)_t^2(x,s)G(x,1)\,dx\bigg)\,ds\\
  &+ \frac12 \int_{\R^N}\bigg(|\nabla
  v_i(x,t)|^2-\frac{a(x/|x|)}{|x|^2}\,v_i^2(x,t)-th(\sqrt t x,
  t+a_i)v_i^2(x,t)\bigg)
  G(x,1)\,dx\\
  &=\frac12 \int_{\R^N}\bigg(|\nabla
  v_{0,i}(x)|^2-\frac{a(x/|x|)}{|x|^2}\,v^2_{0,i}(x)- T_i h(\sqrt {
    T_i} x, T_i+a_i)v_{0,i}^2(x)\bigg) G(x,1)\,dx
  \\
  &+\int_t^{T_i}\!\!\bigg(\int_{\R^N}h(\sqrt s x, s+a_i)\bigg(
  \frac{|x|^2}8v_i^2(x,s)- \frac{\nabla v_i(x,s)\cdot x }{2}v_i(x,s)
  -\frac{N-2}{4}v_i^2(x,s)\bigg)G(x,1)\,dx\bigg)\,ds\\
  &+\frac{1}{2}\int_t^{T_i}s\bigg(\int_{\R^N}h_{s}(\sqrt s x,
  s+a_i)v_i^2(x,s)G(x,1)\,dx\bigg)\,ds,
\end{align*}
for all $t\in (0, T_i)$, where $v_{0,i}(x):=u_i(\sqrt { T_i} x, T_i)
\in{\mathcal H}$.  Therefore the function
$$
t\mapsto \int_{\R^N}\bigg(|\nabla
v_i(x,t)|^2-\frac{a(x/|x|)}{|x|^2}\,v_i^2(x,t)-th(\sqrt t x,
t+a_i)v_i^2(x,t)\bigg) G(x,1)\,dx
$$
is absolutely continuous  in $(\tau,T_i)$ for all $\tau\in(0,T_i)$ and
\begin{align*}
  \frac{d}{dt}&\int_{\R^N}\bigg(|\nabla
  v_i(x,t)|^2-\frac{a(x/|x|)}{|x|^2}\,v_i^2(x,t)-th(\sqrt t x,
  t+a_i)v_i^2(x,t)\bigg)
  G(x,1)\,dx\\
  &=2t \int_{\R^N}(v_i)_t^2(x,t)G(x,1)\,dx\\
  &\quad-\int_{\R^N}h(\sqrt t x, t+a_i)\bigg( \frac{|x|^2}4v_i^2(x,t)-
  (\nabla v_i(x,t)\cdot x )v_i(x,t)
  -\frac{N-2}{2}v_i^2(x,t)\bigg)G(x,1)\,dx\\
  &-t\int_{\R^N}h_{s}(\sqrt s x, s+a_i)v_i^2(x,s)G(x,1)\,dx.
\end{align*}
The change of variables $u_i(x,t)=v_i(x/\sqrt t,t)$ leads to the
conclusion in case {\bf (I)}.

Let us now consider case {\bf (II)}, i.e. $f (x,t, u)=\varphi (x,t,
u)$ with $\varphi$ satisfying \eqref{eq:fi} and $u$ satisfying
(\ref{eq:u1}--\ref{eq:u2}). We test equation (\ref{eq:eqforv_i}) with
$(v_i)_t$ (passing through a suitable approximation) and, by Corollary
\ref{c:pos_per_nonlin}, we obtain, for all $t\in(0,T_i)$,
\begin{align*}
  \int_t^{T_i}s&\bigg(\int_{\R^N}(v_i)_t^2(x,s)G(x,1)\,dx\bigg)\,ds\leq{\rm
    const\,}\bigg( \|u_i(\sqrt {T_i}\,\cdot,T_i)\|^2_{\mathcal H}+
  \int_{\R^N}v_i^2(x,t)G(x,1)\,dx\bigg)\\
  &\quad -\int_t^{T_i}s\bigg(\int_{\R^N}\varphi(\sqrt s x, s+a_i,
  v_i(x,t))(v_i)_{t}(x,t)G(x,1)\,dx\bigg)\,ds\\
  &\quad+\frac{1}{2}\int_t^{T_i}\frac{d}{ds}\bigg(s\int_{\R^N}\varphi(\sqrt
  s x, s+a_i,v_i(x,t))v_i(x,s)G(x,1)\,dx\bigg)\,ds.
\end{align*}
Since in view of hypothesis \eqref{eq:fi} on $\varphi$, conditions
\eqref{eq:u1} and \eqref{eq:u2} on $u$, and Lemma \ref{l:sob} the
integrals at the right hand side lines of the previous formula are finite for
every $t\in(0,T_i)$, we conclude that
$$
(v_i)_t\in L^2(\tau,T_i;{\mathcal L})\quad\text{for all
  }\tau\in (0,T_i).
$$
Testing (\ref{eq:eqforv_i}) for $v_i$ with $(v_i)_t$ also yields
\begin{align*}
  \int_t^{T_i}s&\bigg(\int_{\R^N}(v_i)_t^2(x,s)G(x,1)\,dx\bigg)\,ds\\
  &\quad+ \frac12 \int_{\R^N}\bigg(|\nabla
  v_i(x,t)|^2-\frac{a(x/|x|)}{|x|^2}\,v_i^2(x,t)-t\varphi(\sqrt t x,
  t+a_i, v_i(x,t))v_i(x,t)\bigg)
  G(x,1)\,dx\\
  &=\frac12 \int_{\R^N}\bigg(|\nabla
  v_{0,i}(x)|^2-\frac{a(x/|x|)}{|x|^2}\,v^2_{0,i}(x)- T_i
  \varphi(\sqrt { T_i} x, T+a_i, v_{0,i}) v_{0,i}(x)\bigg) G(x,1)\,dx
  \\
  &\quad -\int_t^{T_i}s\bigg(\int_{\R^N}
  \varphi(\sqrt s x, s+a_i, v_i(x,t))(v_i)_{t}(x,t)G(x,1)\,dx\bigg)\,ds\\
  &\quad+\frac{1}{2}\int_t^{T_i}\frac{d}{ds}\bigg(s\int_{\R^N}\varphi(\sqrt
  s x, s+a_i,v_i(x,s))v_i(x,s)G(x,1)\,dx\bigg)\,ds,
\end{align*}
for a.e. $t\in (0, T_i)$, where $v_{0,i}(x):=u_i(\sqrt { T_i} x, T_i)
\in{\mathcal H}$.  Therefore the function
$$
t\mapsto \int_{\R^N}\bigg(|\nabla
v_i(x,t)|^2-\frac{a(x/|x|)}{|x|^2}\,v_i^2(x,t)-t\varphi(\sqrt t x, t+a_i,
v_i(x,t))v_i(x,t)\bigg) G(x,1)\,dx
$$
is absolutely continuous  in $(0,\tau)$ for all $\tau\in(0,T_i)$ and
\begin{align*}
  \frac{d}{dt}&\int_{\R^N}\bigg(|\nabla
  v_i(x,t)|^2-\frac{a(x/|x|)}{|x|^2}\,v_i^2(x,t)-t\varphi(\sqrt t x,
  t+a_i, v_i(x,t))v_i (x,t)\bigg)
  G(x,1)\,dx\\
  &=2t \int_{\R^N}(v_i)_t^2(x,t)G(x,1)\,dx +2t\int_{\R^N}
  \varphi(\sqrt t x, t+a_i, v_i(x,t))(v_i)_{t}(x,t)G(x,1)\,dx\\
  &\quad-\frac{d}{dt}\bigg(t\int_{\R^N}\varphi(\sqrt t x,
  t+a_i,v_i(x,t))v_i(x,t)G(x,1)\,dx\bigg).
\end{align*}
The change of variables $u_i(x,t)=v_i(x/\sqrt t,t)$ leads to
\begin{gather*}
  \frac{d}{dt}(t D_i(t))=2t
  \int_{\R^N}\left|(u_i)_t(x,t)+\frac{\nabla u_i(x,t)\cdot
      x}{2t}\right|^2G(x,t)\,dx\\
   +2t\int_{\R^N} \varphi(x, t+a_i, u_i(x,t))
  (u_i)_t(x,t)G(x,t)\,dx+\int_{\R^N} \varphi(x, t+a_i, u_i(x,t))
  \nabla u_i(x,t)\cdot x
  G(x,t)\,dx\\
   -\frac{d}{dt}\bigg(t\int_{\R^N}\varphi( x,
  t+a_i,u_i(x,t))u_i(x,t)G(x,t)\,dx\bigg)
\end{gather*}
and hence
\begin{gather*}
  \frac{d}{dt}(t D_i(t))=2t
  \int_{\R^N}\left|(u_i)_t(x,t)+\frac{\nabla u_i(x,t)\cdot
      x}{2t}\right|^2G(x,t)\,dx\\
   +2t\int_{\R^N} \varphi(x, t+a_i, u_i(x,t))
  (u_i)_t(x,t)G(x,t)\,dx+\int_{\R^N} \varphi(x, t+a_i, u_i(x,t))
  \nabla u_i(x,t)\cdot x\,
  G(x,t)\,dx\\
+\frac{N-2}2\int_{\R^N}\varphi( x,
  t+a_i,u_i(x,t))u_i(x,t)G(x,t)\,dx
-t \int_{\R^N}\frac{\partial \varphi}{\partial t}( x,
  t+a_i,u_i(x,t))u_i(x,t)G(x,t)\,dx\\
-t \int_{\R^N}\bigg(\frac{\partial \varphi}{\partial u_i}( x,
  t+a_i,u_i(x,t))u_i(x,t)+\varphi( x,
  t+a_i,u_i(x,t))\bigg)
(u_i)_t(x,t)G(x,t)\,dx\\
-\int_{\R^N}\frac{|x|^2}{4t}\varphi( x,
  t+a_i,u_i(x,t))u_i(x,t)G(x,t)\,dx.
\end{gather*}
Integration by parts yields (these formal computations can be made
rigorous through a suitable approximation)
\begin{gather*}
  \int_{\R^N} \varphi(x, t+a_i, u_i(x,t)) \nabla u_i(x,t)\cdot x\,
  G(x,t)\,dx=-N\int_{\R^N} \Phi(x,t+a_i, u_i(x,t))
  G(x,t)\,dx\\
  + \int_{\R^N} \frac{|x|^2}{2t}\Phi(x,t+a_i, u_i(x,t)) G(x,t)\,dx-
  \int_{\R^N} \nabla_x\Phi(x,t+a_i, u_i(x,t))\cdot x G(x,t)\,dx
\end{gather*}
thus yielding the conclusion in case {\bf (II)}.
\end{pf}

For all $i=1,\dots,k$, let us introduce the \emph{Almgren type
  frequency function} associated to $u_i$
\begin{equation}\label{eq:9}
N_i:(0,2\alpha)\to\R\cup\{-\infty,+\infty\},
\quad N_i(t):=\frac{tD_i(t)}{H_i(t)}.
\end{equation}
Frequency functions associated to unperturbed parabolic equations of
type (\ref{prob}) (i.e. in the case $f(x,t,s)\equiv 0$) were first
studied by C.-C. Poon in \cite{poon}, where unique continuation
properties are derived by proving monotonicity of the quotient in
(\ref{eq:9}).  Due to the presence of the perturbing function
$f(x,t+a_i,u(x,t))$, the functions $N_i$ will not be nondecreasing as
in the case treated by Poon; however in both cases {\bf (I)} and {\bf
  (II)}, we can prove that their derivatives are integrable
perturbations of nonnegative functions wherever the $N_i$'s assume
finite values.  Moreover our analysis will show that actually the
$N_i$'s assume finite values all over $(0,2\alpha)$.

\begin{Lemma}\label{l:Nprime}
  Let $i\in \{1,\dots,k\}$. If there exist $\beta_i,T_i\in
  (0,2\alpha)$ such that
\begin{equation}\label{eq:10}
\beta_i<T_i, \quad H_i(t)>0 \text{ for all $t\in (\beta_i,T_i)$},\quad
\text{and}\quad u_i(\cdot,T_i)\in {\mathcal H}_{T_i},
\end{equation}
then the function $N_i$ defined in \eqref{eq:9} belongs to
$W^{1,1}_{\rm loc}(\beta_i,T_i)$ and
\begin{align*}
N'_i(t)={\nu}_{1i}(t)+{\nu}_{2i}(t)
\end{align*}
in a distributional sense and a.e. in $(\beta_i,T_i)$ where
\begin{align*}
  {\nu}_{1i}(t)&=\frac{2t}{H_{i}^2(t)}{{
      \Bigg(\bigg(\int_{\R^N}\bigg|(u_i)_t(x,t)+\frac{\nabla u_i(x,t)\cdot
        x}{2t}\bigg|^2G(x,t)\,dx\bigg)
      \bigg(\int_{\R^N}u_{i}^2(x,t)\, G(x,t)\,dx\bigg)}}\\
  &\hskip4cm{{-\bigg(\int_{\R^N}\Big((u_i)_t(x,t)+\frac{\nabla
        u_{i}(x,t)\cdot x}{2t}\Big)u_{i}(x,t)G(x,t)\,dx
      \bigg)^{\!2}\Bigg)}}
\end{align*}
and ${\nu}_{2i}$ is as follows:

\smallskip\noindent
{in  case \bf (I)}
\begin{align*} {\nu}_{2i}(t)&={{\dfrac1{H_{i}(t)}
      \int_{\R^N}h(x,t+a_i)\left(\frac{N-2}{2}u_{i}^2(x,t)+(\nabla
        u_{i}(x,t)\cdot x)u_i(x,t)-\frac{|x|^2}{4t}u_{i}^2(x,t)\right)
      G(x,t)\,dx}}\\
  &\quad{{-\dfrac
      t{H_{i}(t)}\bigg(\int_{\R^N}h_{t}(x,t+a_i)u_{i}^2(x,t)G(x,t)\,dx\bigg)}},
\end{align*}
{in case \bf (II)}
\begin{align*}
  {\nu}_{2i}(t)=\,&{\dfrac1{H_{i}(t)}} \bigg(
t\int_{\R^N} \Big(
  \varphi(x,t+a_i,u_i(x,t))-\frac{\partial\varphi}{\partial u_i}
  (x,t+a_i,u_i(x,t))u_i(x,t)\Big)(u_i)_t(x,t)G(x,t)\,dx\\
  &\hskip1cm+\int_{\R^N} \Big(\frac{N-2}2\varphi(x,t+a_i,u_i(x,t))u_i(x,t)-t
  \frac{\partial\varphi}{\partial t}(x,t+a_i,u_i(x,t))
  u_i(x,t)\\
&\hskip3cm-N\Phi(x,t+a_i,u_i(x,t))-\nabla_x\Phi(x,t+a_i,u_i(x,t))\cdot
  x\Big)G(x,t)\,dx\\
  &\hskip1cm+\int_{\R^N}\frac{|x|^2}{4t} \Big(2\Phi(x,t+a_i,u_i(x,t))-
  \varphi(x,t+a_i,u_i(x,t))u_i(x,t)\Big)G(x,t)\,dx\bigg).
\end{align*}
\end{Lemma}
\begin{pf}
  From Lemma \ref{l:Hprime} and \ref{l:Dprime}, it follows that
  $N_{i}\in W^{1,1}_{\rm loc}(\beta_i,T_i)$. From (\ref{eq:10i}) we
  deduce that
$$
N'_{i}(t)=\frac{(tD_{i}(t))'H_{i}(t)-tD_{i}(t)H'_{i}(t)}{H_{i}^2(t)}=
\frac{(tD_{i}(t))'H_{i}(t)-2tD_{i}^2(t)}{H_{i}^2(t)},
$$
which yields the conclusion in view of (\ref{eq:Hi(t)}),
(\ref{eq:Di(t)}), and Lemma \ref{l:Dprime}.
\end{pf}

\noindent The term ${\nu_{2i}}$ can estimated as follows.
\begin{Lemma}\label{l:est_N_2}
There exists $C_3>0$ such that, if $i\in \{1,\dots,k\}$ and $\beta_i,T_i\in
  (0,2\alpha)$ satisfy (\ref{eq:10}), then
\begin{equation*}
  \Big|{\nu_{2i}}(t)\Big|\leq
\begin{cases}
  C_3\big(N_i(t)+{\textstyle{\frac{N-2}{4}}}\big)
\big(t^{-1+\e/2}+\|h_t(\cdot,t+a_i)\|_{L^{N/2}(\R^N)}\big),
  &\text{in case {\bf (I)}},\\[5pt]
  C_3\big(N_i(t)+{\textstyle{\frac{N-2}{4}}}\big)
\,t^{-1+\frac {N+2-p(N-2)}{2(p+1)}}, &\text{in case
    {\bf (II)} if $i=1$},\\[5pt]
  C_3\beta_i^{-1}\big(N_i(t)+{\textstyle{\frac{N-2}{4}}}\big)
\,t^{-1+\frac {N+2-p(N-2)}{2(p+1)}}, &\text{in case
    {\bf (II)} if $i>1$},
\end{cases}
\end{equation*}
for a.e. $t\in(\beta_i,T_i)$, where ${\nu_{2i}}$ is as in Lemma \ref{l:Nprime}.
\end{Lemma}
\begin{pf}
  Let us first consider case {\bf (I)}, i.e. $f (x,t, u)=h(x,t)u$,
  with $h(x,t)$ under conditions (\ref{eq:der}--\ref{eq:h}).  In order
  to estimate ${\nu}_{2i}$ we observe that, from (\ref{eq:h}),
\begin{align}\label{eq:15}
  &\bigg|\int_{\R^N}h(x,t+a_i)(\nabla u_{i}(x,t)\cdot
  x)u_{i}(x,t) G(x,t)\,dx\bigg|\\
  &\notag\leq
  C_h\int_{\R^N}(1+|x|^{-2+\e})|\nabla u_{i}(x,t)||x||u_{i}(x,t) |G(x,t)\,dx\\
  \notag &\leq C_ht \int_{\R^N}|\nabla
  u_{i}(x,t)|\frac{|x|}{t}|u_{i}(x,t)| G(x,t)\,dx +C_ht^{\e/2}
  \int_{\{|x|\leq \sqrt t\}}|\nabla u_{i}(x,t)|\frac{|u_{i}(x,t)|}{|x|} G(x,t)\,dx\\
  \notag &\hskip4cm +C_ht^{\e/2}\int_{\{|x|\geq \sqrt t\}}|\nabla
  u_{i}(x,t)|\frac{|x|}{t}|u_{i}(x,t)| G(x,t)\,dx
  \\
  \notag & \leq \frac12C_h(t+t^{\e/2})\int_{\R^N}|\nabla
  u_{i}(x,t)|^2G(x,t)\,dx+
  \frac12C_h(t+t^{\e/2})\int_{\R^N}\frac{|x|^2}{t^2}u_{i}^2(x,t)G(x,t)\,dx\\
  \notag &\hskip2cm + \frac12C_ht^{\e/2}\int_{\R^N}|\nabla
  u_{i}(x,t)|^2G(x,t)\,dx
  + \frac12C_ht^{\e/2}\int_{\R^N}\frac{u_{i}^2(x,t)}{|x|^2}G(x,t)\,dx\\
  \notag&\leq \frac12C_ht^{\e/2}(2+{\overline T}^{1-\e/2})
  \int_{\R^N}|\nabla u_{i}(x,t)|^2G(x,t)\,dx\\
  \notag &\quad +
  \frac12C_ht^{\e/2}(1+{\overline T}^{1-\e/2})\int_{\R^N}
  \frac{|x|^2}{t^2}u_{i}^2(x,t)G(x,t)\,dx +
  \frac12C_ht^{\e/2}\int_{\R^N}\frac{u_{i}^2(x,t)}{|x|^2}G(x,t)\,dx ,
\end{align}
and
\begin{align}\label{eq:14}
&\int_{\R^N}|h(x,t+a_i)||x|^2u_{i}^2(x,t) G(x,t)\,dx\leq C_h
\int_{\R^N}|x|^2u_{i}^2(x,t) G(x,t)\,dx\\
\notag&\hskip7cm+C_h
\int_{\R^N}|x|^{-2+\e}|x|^2u_{i}^2(x,t) G(x,t)\,dx\\
\notag &\leq  C_h
\int_{\R^N}|x|^2u_{i}^2(x,t) G(x,t)\,dx+C_h t^{\e/2}
\int_{\{|x|\leq \sqrt t\}}u_{i}^2(x,t) G(x,t)\,dx\\
\notag &
\hskip7cm+C_ht^{-1+\e/2}
\int_{\{|x|\geq \sqrt t \}}|x|^2u_{i}^2(x,t) G(x,t)\,dx\\
\notag &\leq  C_ht^{-1+\e/2}(1+{\overline T}^{1-\e/2})
\int_{\R^N}|x|^2u_{i}^2(x,t) G(x,t)\,dx+C_h t^{\e/2}
\int_{\R^N}u_{i}^2(x,t) G(x,t)\,dx,
\end{align}
for a.e. $t\in (\beta_i,T_i)$. Moreover, by H\"older's inequality and
Corollary \ref{cor:ineqSob},
\begin{align}\label{eq:hip}
&\int_{\R^N}|h_t(x,t+a_i)|u_{i}^2(x,t) G(x,t)\,dx\leq C_{2^*}t^{-1}
\|u_i\|^2_{{\mathcal H}_t}\|h_t(\cdot,t+a_i)\|_{L^{{N}/{2}}(\R^N)}
\end{align}
for a.e. $t\in (\beta_i,T_i)$.
Collecting (\ref{eq:11}), (\ref{eq:15}),
(\ref{eq:14}) and (\ref{eq:hip}), we obtain that
\begin{multline}\label{eq:16}
\Big|{\nu_{2i}}(t)\Big|\leq \frac{{\rm const}\, t^{\e/2}}{H_{i}(t)}
\bigg(\frac1t\int_{\R^N}u_{i}^2(x,t) G(x,t)\,dx+
\int_{\R^N}\frac{u_{i}^2(x,t)}{|x|^2}G(x,t)\,dx\\
+\int_{\R^N}|\nabla u_{i}(x,t)|^2G(x,t)\,dx+
\frac1{t^2}\int_{\R^N}|x|^2u_{i}^2(x,t) G(x,t)\,dx\bigg)\\
+\frac{C_{2^*}}{H_{i}(t)}
\|u_i\|^2_{{\mathcal H}_t}\|h_t(\cdot,t+a_i)\|_{L^{{N}/{2}}(\R^N)}.
\end{multline}
From inequality (\ref{eq:16}), Lemma \ref{Hardytemp}, Corollary
\ref{c:pos_per}, and Corollary \ref{cor:ineq}, we deduce that there
exists $C_3>0$ depending only on $C_h$, ${\overline T}$, and $N$, such
that, for a.e. $t\in(\beta_i,T_i)$,
\begin{align*}
\Big|{\nu_{2i}}(t)\Big|&\leq
\frac{C_3}{H_{i}(t)}\Big(tD_i(t)+{\textstyle{\frac{N-2}{4}}}H_i(t)\Big)\Big(
t^{-1+\e/2}+\|h_t(\cdot,t+a_i)\|_{L^{{N}/{2}}(\R^N)}\Big)\\
&=C_3\big(N_i(t)+{\textstyle{\frac{N-2}{4}}}\big)\Big(
t^{-1+\e/2}+\|h_t(\cdot,t+a_i)\|_{L^{{N}/{2}}(\R^N)}\Big)
\end{align*}
thus completing the proof in case {\bf (I)}.

Let us now consider case {\bf (II)}, i.e. $f (x,t, s)=\varphi (x,t,
s)$ with $\varphi$ under condition \eqref{eq:fi} and $u$ satisfying
\eqref{eq:u1} and \eqref{eq:u2}. From (\ref{eq:fi}), we have that
\begin{gather}\label{eq:84}
  \Big|{\nu_{2i}}(t)\Big|\leq\frac{\rm const}{H_i(t)}
  \bigg(t\int_{\R^N}|u_i(x,t)|^q|(u_i)_t(x,t)|G(x,t)\,dx\\
  \notag\quad+\int_{\R^N}\big(|u_i(x,t)|^2+|u_i(x,t)|^{p+1}\big)G(x,t)\,dx
  +\int_{\R^N}\frac{|x|^2}{t}\big(|u_i(x,t)|^2+|u_i(x,t)|^{p+1}\big)G(x,t)\,dx
  \bigg).
\end{gather}
From H\"older's inequality, Corollary \ref{cor:ineqSob}, and assumptions
(\ref{eq:u1}--\ref{eq:u2}), it follows that
\begin{align}\label{eq:85}
t&\int_{\R^N}|u_i(x,t)|^q|(u_i)_t(x,t)|G(x,t)\,dx\\
\notag&\leq t \bigg(\int_{\R^N}|u_i(x,t)|^{p+1}G^{\frac{p+1}{2}}(x,t)\,dx
  \bigg)^{\!\!\frac{2}{p+1}}\|u(\cdot,t+a_i)\|^{q-2}_{L^{p+1}(\R^N)}
\|u_t(\cdot,t+a_i)\|_{L^{\frac{p+1}{p+1-q}}(\R^N)}\\
\notag&\leq {\rm const\,}t^{-\frac N{p+1}\frac{p-1}2}\|u_i\|^2_{{\mathcal H}_t}
\end{align}
and, taking into account also Corollary \ref{cor:ineq},
\begin{align}\label{eq:86}
  \int_{\R^N}&\frac{|x|^2}{t}\big(|u_i(x,t)|^2+|u_i(x,t)|^{p+1}\big)G(x,t)\,dx
  \leq \int_{\R^N}\frac{|x|^2}{t}|u_i(x,t)|^2 G(x,t)\,dx
  \\
\notag  &+
\frac{t+a_i}{t}
\bigg(\int_{\R^N}|u_i(x,t)|^{p+1}G^{\frac{p+1}{2}}(x,t)\,dx
  \bigg)^{\!\!\frac{2}{p+1}} \bigg(\int_{\R^N}\bigg(\frac{|x|^2}{t+a_i}
  \bigg)^{\!\!\frac{p+1}{p-1}}|u(x,t+a_i)|^{p+1}\,dx
  \bigg)^{\frac{p-1}{p+1}}\\[8pt]
  \notag&\hskip3cm\leq
  \begin{cases} {\rm const\,}t^{-\frac
      N{p+1}\frac{p-1}2}\|u_i\|^2_{{\mathcal
        H}_t},&\text{if }i=1,\\[8pt]
    {\rm const\,}b_i\beta_i^{-1}t^{-\frac
      N{p+1}\frac{p-1}2}\|u_i\|^2_{{\mathcal H}_t},&\text{if }i>1.
  \end{cases}
\end{align}
As in (\ref{eq:varphi}) we can estimate
\begin{align}\label{eq:87}
\int_{\R^N}\big(|u_i(x,t)|^2+|u_i(x,t)|^{p+1}\big)G(x,t)\,dx
\leq  {\rm const\,}t^{-\frac N{p+1}\frac{p-1}2}\|u_i\|^2_{{\mathcal H}_t}.
\end{align}
Collecting (\ref{eq:84}), (\ref{eq:85}), (\ref{eq:86}), and (\ref{eq:87}),
and using Corollary \ref{c:pos_per_nonlin}, we obtain that
there
exists some positive constant $C_3$ such
that, for a.e. $t\in(\beta_i,T_i)$,
\begin{align*}
  \Big|{\nu_{2i}}(t)\Big|\leq
  \begin{cases}
    \frac{C_3}{H_{i}(t)}\,t^{-\frac N{p+1}\frac{p-1}2}\big(tD_i(t)
    +{\textstyle{\frac{N-2}{4}}}H_i(t)\big)=C_3\big(N_i(t)+
    {\textstyle{\frac{N-2}{4}}}\big)\,t^{-1+\frac
      {N+2-p(N-2)}{2(p+1)}},&\!\!\!\text{if }i=1,\\[8pt]
    \frac{C_3\beta_i^{-1}}{H_{i}(t)}\,t^{-\frac
      N{p+1}\frac{p-1}2}\big(tD_i(t)+{\textstyle{\frac{N-2}{4}}}
    H_i(t)\big)=C_3\beta_i^{-1}
    \big(N_i(t)+{\textstyle{\frac{N-2}{4}}}\big)\,t^{-1+\frac
      {N+2-p(N-2)}{2(p+1)}}\!\!,&\!\!\!\text{if }i>1,
\end{cases}
\end{align*}
thus completing the proof in case {\bf (II)}.
\end{pf}
\begin{Lemma}\label{l:Nabove}
There exists $C_4>0$ such that, if $i\in \{1,\dots,k\}$ and $\beta_i,T_i\in
  (0,2\alpha)$ satisfy (\ref{eq:10}), then,
 for every $t\in(\beta_i,T_i)$,
\begin{equation*}
  N_i(t)\leq
\begin{cases}
-\frac{N-2}{4}+C_{4}\big(N_i(T_i)+\frac{N-2}4\big),
&\text{in case {\bf (I)} and in case {\bf (II)} if $i=1$},\\[8pt]
-\frac{N-2}{4}+C_{4}^{1/\beta_i}\big(N_i(T_i)+\frac{N-2}4\big),
&\text{in case  {\bf (II)} if $i>1$}.
\end{cases}
\end{equation*}
\end{Lemma}
\begin{pf}
Let ${\nu}_{1i}$ and ${\nu}_{2i}$ as in Lemma \ref{l:Nprime}.
By Schwarz's inequality,
\begin{equation}\label{eq:17}
{\nu}_{1i}\geq 0 \quad\text{a.e. in }(\beta_i,T_i).
\end{equation}
From Lemma \ref{l:Nprime}, (\ref{eq:17}), and Lemma \ref{l:est_N_2},
 we deduce that
$$
\frac{d}{dt}N_{i}(t)\geq
\begin{cases}
 - C_3\big(N_i(t)+{\textstyle{\frac{N-2}{4}}}\big)
\Big(t^{-1+\e/2}+\|h_t(\cdot,t+a_i)\|_{L^{N/2}(\R^N)}\Big),
  &\text{in case {\bf (I)}},\\[8pt]
  -C_3\big(N_i(t)+{\textstyle{\frac{N-2}{4}}}\big)
\,t^{-1+\frac {N+2-p(N-2)}{2(p+1)}}, &\text{in case
    {\bf (II)} if $i=1$},\\[8pt]
  -C_3\beta_i^{-1}\big(N_i(t)+{\textstyle{\frac{N-2}{4}}}\big)
\,t^{-1+\frac {N+2-p(N-2)}{2(p+1)}}, &\text{in case
    {\bf (II)} if $i>1$},
\end{cases}
$$
for a.e. $t\in(\beta_i,T_i)$.  After integration, it follows that
\begin{multline*}
  N_{i}(t)\\
  \leq
\begin{cases}
-\frac{N-2}{4}
+\Big(
N_{i}(T_{i})+\frac{N-2}{4}\Big)\exp\Big({\frac{2C_3}{\e}T_{i}^{\e/2}+C_3
\|h_t\|_{L^{1}((0,T),L^{{N}/{2}}(\R^N))}}\Big),
  &\text{in case {\bf (I)}},\\[8pt]
-\frac{N-2}{4}
+\Big(
N_{i}(T_{i})+\frac{N-2}{4}\Big)
\exp\Big(\frac{2(p+1)C_3}{N+2-p(N-2)}T_{i}^{\frac{N+2-p(N-2)}{2(p+1)}}\Big),
  &\text{in case {\bf (II)}, $i=1$},\\[8pt]
-\frac{N-2}{4}
+\Big(
N_{i}(T_{i})+\frac{N-2}{4}\Big)
\exp\Big(\frac{2(p+1)C_3\beta_i^{-1}}{N+2-p(N-2)}T_{i}^{\frac{N+2-p(N-2)}{2(p+1)}}\Big),
  &\text{in case {\bf (II)}, $i>1$},
\end{cases}
\end{multline*}
for any $t\in (\beta_i,T_i)$, thus yielding the conclusion.
\end{pf}

\begin{Lemma}\label{l:t_0=0}
Let $i\in\{1,\dots,k\}$. If $H_i\not\equiv 0$, then
$$
H_i(t)>0
\quad\text{for all }t\in(0,2\alpha).
$$
\end{Lemma}
\begin{pf}
  From continuity of $H_i$, the assumption $H_i\not\equiv0$, and the
  fact that $u_i(\cdot,t)\in {\mathcal H}_t$ for a.e. $t\in (0,2\alpha)$,
we deduce that there exists $T_i\in (0,2\alpha)$ such that
\begin{equation}\label{eq:78}
H_i(T_i)>0\quad\text{and} \quad u_i(\cdot,T_i)\in {\mathcal H}_{T_i}.
\end{equation}
Lemma \ref{l:Hpos}
implies that $H_i(t)>0$ for all $t\in [T_i,2\alpha)$.
 We consider
 \begin{equation*}
t_i:=\inf\{s\in (0,T_i):H_i(t)>0\text{ for all }t\in(s,2\alpha)\}.
\end{equation*}
Due to  Lemma \ref{l:Hpos},
either
\begin{equation}\label{eq:40}
t_i=0 \text{ and } H_i(t)>0\text{ for all }t\in (0,2\alpha)
\end{equation}
or
\begin{equation}\label{eq:39}
0<t_i<T_i\text{ and }
\begin{cases}
H_i(t)=0&\text{if }t\in(0,t_i]\\
H_i(t)>0&\text{if }t\in(t_i,2\alpha)
\end{cases}.
\end{equation}
The argument below will exclude alternative \eqref{eq:39}.
Assume by contradiction that \eqref{eq:39} holds. 
From Lemma
\ref{l:Nabove} and  (\ref{eq:10i}), it follows
$$
\frac t2H'_{i}(t)\leq  c_i H_{i}(t)
$$
where
$$
c_i=
\begin{cases}
-\frac{N-2}{4}+C_{4}\big(N_i(T_i)+\frac{N-2}4\big)
&\text{in case {\bf (I)} and in case {\bf (II)} if $i=1$},\\[8pt]
-\frac{N-2}{4}+C_{4}^{1/t_i}\big(N_i(T_i)+\frac{N-2}4\big)
&\text{in case  {\bf (II)} if $i>1$},
\end{cases}
$$
for a.e. $t\in (t_i,T_i)$. By integration, it follows that
\begin{equation}\label{eq:18}
  H_{i}(t)\geq \frac{H_{i}(T_i)}{T_i^{2c_i}}\,t^{2c_i}
  \quad\text{for all
  }t\in[t_i,T_i).
\end{equation}
By (\ref{eq:39}) $H_{i}(t_i)=0$, giving rise to contradiction with
\eqref{eq:18} because of (\ref{eq:78}). Therefore, we exclude
(\ref{eq:39}) and conclude that (\ref{eq:40}) holds.
\end{pf}
\begin{Lemma}\label{l:H_{i+1}}
Let $i\in\{1,\dots,k\}$. Then
$$
H_i(t)\equiv 0\text { in $(0,2\alpha)$ if and only if } H_{i+1}(t)\equiv 0
\text { in $(0,2\alpha)$}.
$$
\end{Lemma}
\begin{pf} First, we prove that $H_i(t)\equiv 0$ in $(0,2\alpha)\text
  { implies } H_{i+1}(t)\equiv 0$ in $(0,2\alpha).$ Let's suppose by
  contradiction that $H_{i+1}(t)\not\equiv 0$. By Lemma \ref{l:t_0=0},
  we conclude that $H_{i+1}(t)>0$ for all $t\in(0,2\alpha)$.  It
  follows that $u_{i+1}(\cdot, t)\not\equiv 0$ for all
  $t\in(0,2\alpha)$ and $u(\cdot, t)\not\equiv 0$, for all
  $t\in(i\alpha,(i+1)\alpha)$. Hence, $u_{i}(\cdot, t)\not\equiv 0$,
  for all $t\in (\alpha, 2\alpha)$ and thus $H_{i}\not\equiv 0$ in
  $(0,2\alpha)$, a contradiction.

  Let us now prove that $H_{i+1}(t)\equiv 0$ in $(0,2\alpha)\text {
    implies } H_{i}(t)\equiv 0$ in $(0,2\alpha).$ Let's suppose by
  contradiction that $H_{i}(t)\not\equiv 0$, then, by Lemma
  \ref{l:Hpos}, $H_{i}(t)>0$ in $(\overline{t},2\alpha)$ for some
  $\overline{t}\in (\alpha,2\alpha)$. Hence, $u_{i}(\cdot,
  t)\not\equiv 0$ in $(\overline{t},2\alpha)$ and then $u_{i+1}(\cdot,
  t)\not\equiv 0$ in $(\overline{t}-\alpha,\alpha)$, thus implying
  $H_{i+1}(t)\not\equiv 0$, a contradiction.
\end{pf}

\begin{Corollary}\label{c:ucontinuation}
  If $u\not\equiv0$ in  $\R^N\times(0,T)$, then
$$
H_i(t)>0
$$
for all $t\in(0,2\alpha)$ and $i=1,\dots,k$. In particular,
\begin{equation}\label{eq:88}
\int_{\R^N}u^2(x,t)G(x,t)\,dx>0\quad\text{for all }t\in(0,T).
\end{equation}
\end{Corollary}
\begin{pf}
  If $u\not\equiv0$, then there exists some $i_0\in\{1,\dots,k\}$
  such that $u_{i_0}\not\equiv 0$ in $(0,2\alpha)$. Hence,
  $H_{i_0}(t)\not\equiv 0$ in $(0,2\alpha)$ and, thanks to lemma \ref
  {l:H_{i+1}}, $H_{i}(t)\not\equiv 0$ in $(0,2\alpha)$ for all
  $i=1,\dots,k$. Applying Lemma \ref{l:t_0=0}, we conclude that, for
  all $i=1,\dots,k$, $H_{i}(t)>0$ in $(0,2\alpha)$, thus implying
  (\ref{eq:88}).
\end{pf}

\begin{pfn}{Proposition \ref{p:uniq_cont}}
It follows immediately from Corollary \ref{c:ucontinuation}.
\end{pfn}

\noindent
Henceforward, we assume $u\not\equiv 0$ and we denote, for all
$t\in(0,2\alpha)$,
\begin{align*}
&H(t)=H_{1}(t)=\int_{\R^N}u^2(x,t)\, G(x,t)\,dx,\\
& D(t)=D_{1}(t)=\!\!\int_{\R^N}\!\!\bigg(|\nabla u(x,t)|^2-
  \dfrac{a\big(\frac{x}{|x|}\big)}{|x|^2}u^2(x,t)-
f(x,t,u(x,t))u(x,t)\bigg)G(x,t)\,dx.
\end{align*}
Corollary \ref{c:ucontinuation} ensures that, if $u\not\equiv 0$ in
$\R^N\times (0,T)$, $H(t)>0$ for all $t\in(0,2\alpha)$ and hence the
\emph{Almgren type frequency function}
$$
{\mathcal N}(t)={\mathcal N}_{f,u}(t)=N_1(t)=\frac{tD(t)}{H(t)}
$$
is well defined over all $(0,2\alpha)$. Moreover, by Lemma \ref{l:Nprime},
${\mathcal N}\in W^{1,1}_{\rm loc}(0,2\alpha)$ and
$$
{\mathcal N}'(t)={\nu}_1(t)+{\nu}_2(t) \quad\text{for a.e. }t\in(0,2\alpha),
$$
where
\begin{equation}\label{eq:nu1nu2}
{\nu}_1(t)={\nu}_{11}(t) \quad\text{and}\quad {\nu}_2(t)={\nu}_{21}(t),
\end{equation}
with
${\nu}_{11},{\nu}_{21}$ as in Lemma \ref{l:Nprime}.
 Since, by (\ref{eq:defsol1}),
 $u(\cdot,t)\in
{\mathcal H}_t$ for a.e. $t\in (0,T)$, we can fix $T_0$ such that
\begin{equation}\label{eq:T_0}
T_0\in (0,2\alpha)
\quad\text{and}\quad
u(\cdot,T_0)\in
{\mathcal H}_{T_0}.
\end{equation}
The following
result clarifies the behavior of ${\mathcal N}(t)$ as $t\to 0^+$.

\begin{Lemma}\label{l:limit}
The limit
$$
\gamma:=\lim_{t\to 0^+}{\mathcal N}(t)
$$
exists and it is finite.
\end{Lemma}
\begin{pf}
  We first observe that ${\mathcal N}(t)$ is bounded from below in
  $(0,2\alpha)$. Indeed
  from Corollaries \ref{c:pos_per} and \ref{c:pos_per_nonlin}, we
  obtain that, for all $t\in(0,2\alpha)$,
$$
tD(t)\geq \bigg(C_1-\frac{N-2}4\bigg)H(t),
$$
and hence
\begin{equation}\label{cota}
{\mathcal N}(t)\geq  C_1-\frac{N-2}4.
\end{equation}
Let $T_0$ as in (\ref{eq:T_0}).  By Schwarz's inequality,
${\nu}_1(t)\geq 0$ for a.e. $t\in (0,T_0)$. Furthermore, from Lemmas
\ref{l:est_N_2} and \ref{l:Nabove}, ${\nu}_2$ belongs to
$L^{1}(0,T_0)$. In particular, ${\mathcal N}'(t)$ turns out to be the
sum of a nonnegative function and of a $L^1$ function over $(0,T_0)$.
Therefore,
$$
{\mathcal N}(t)={\mathcal N}(T_0)-\dint_{t}^{T_0}{\mathcal N}'(s)\,ds
$$
admits a limit as $t\rightarrow 0^{+}$ which is finite in view of
\eqref{cota} and Lemma \ref{l:Nabove}.
\end{pf}

\begin{Lemma}\label{stimaH}
  Let $\gamma:=\lim_{t\rightarrow 0^+} {\mathcal N}(t)$ be as in Lemma
  \ref{l:limit}.  Then there exists a constant $K_1>0$ such that
\begin{equation}\label{eq:52}
H(t)\leq K_1 t^{2\gamma}  \quad \text{for all } t\in (0,T_0).
\end{equation}
Furthermore, for any $\sigma>0$, there exists a constant
$K_2(\sigma)>0$ depending on $\sigma$ such that
\begin{equation}\label{eq:53}
  H(t)\geq K_2(\sigma)\, t^{2\gamma+\sigma}\quad \text{for all } t\in (0,T_0).
\end{equation}
\end{Lemma}
\begin{pf}
 From Lemma \ref{l:Nprime}, \eqref{eq:17}, Lemma \ref{l:est_N_2}, and
Lemma \ref{l:Nabove},
we infer that
\begin{align*}
  {\mathcal N}(t)-\gamma&=\int_0^t ({\nu}_1(s)+{\nu}_2(s))ds \geq
  \int_0^t {\nu}_2(s)ds\\[5pt]
&\geq
\begin{cases}
-  C_3C_4\big({\mathcal N}(T_0)+{\textstyle{\frac{N-2}{4}}}\big)
\int_0^t
\Big(s^{-1+\e/2}+\|h_t(\cdot,s)\|_{L^{N/2}(\R^N)}\Big)\,ds,
  &\text{in case {\bf (I)}},\\[5pt]
-  C_3C_4\big({\mathcal N}(T_0)+{\textstyle{\frac{N-2}{4}}}\big)
\int_0^t s^{-1+\frac {N+2-p(N-2)}{2(p+1)}}\,ds, &\text{in case
    {\bf (II)}},
\end{cases}\\[5pt]
&\geq
\begin{cases}
-  C_3C_4\big({\mathcal N}(T_0)+{\textstyle{\frac{N-2}{4}}}\big)
\Big(\frac2\e t^{\e/2}+\|h_t\|_{L^{r}((0,T),L^{{N}/{2}}(\R^N))}t^{1-1/r}\Big),
  &\text{in case {\bf (I)}},\\[5pt]
-   C_3C_4\big({\mathcal N}(T_0)+{\textstyle{\frac{N-2}{4}}}\big)\frac{2(p+1)}{N+2-p(N-2)} t^{\frac {N+2-p(N-2)}{2(p+1)}}, &\text{in case
    {\bf (II)}}
\end{cases}
\\[5pt]
&\geq - C_5t^\delta
\end{align*}
with
\begin{equation}\label{eq:delta}
\delta=
\begin{cases}
  \min\{\e/2,1-1/r\},
  &\text{in case {\bf (I)}},\\[5pt]
  \frac {N+2-p(N-2)}{2(p+1)}, &\text{in case {\bf (II)}},
\end{cases}
\end{equation}
for some constant $C_5>0$ and for all $t\in(0,T_0)$. From above and
\eqref{eq:10i}, we deduce that
$$
( \log H(t))'=\frac{H'(t)}{H(t)}=\frac{2}{t}{\mathcal N}(t)\geq
\frac2t\gamma-2 C_5t^{-1+\delta}.
$$
Integrating over $(t, T_0)$ we obtain
$$
H(t)\leq \frac{H(T_0)}{T_0^{2\gamma}}e^{2 C_5T_0^{\delta}}t^{2\gamma}
$$
for all $t\in(0,T_0)$, thus proving \eqref{eq:52}.

Let us prove \eqref{eq:53}. Since $\gamma=\lim_{t\rightarrow 0^+}
{\mathcal N}(t)$, for any $\sigma>0$ there exists $t_\sigma>0$ such
that ${\mathcal N}(t)<\gamma+\sigma/2$ for any $t\in (0,t_\sigma)$ and
hence
$$
\frac{H'(t)}{H(t)}=\frac{2\,{\mathcal N}(t)}{t}<\frac{2\gamma+\sigma}{t}.
$$
Integrating over the interval $(t,t_\sigma)$ and by continuity of $H$
outside $0$, we obtain \eqref{eq:53}  for some constant $K_2(\sigma)$
depending on $\sigma$.
\end{pf}

\section{The blow-up analysis}\label{sec:blow-up-analysis}

If $u$ is a weak solution to (\ref{prob}) in the
sense of Definition \ref{def:solution}, then, for every $\lambda>0$,
the function
\begin{equation*}
u_{\lambda}(x,t)=u(\lambda x,\lambda^2t)
\end{equation*}
is a weak solution to
\begin{equation}\label{lambda}
(u_{\lambda})_t+\Delta u_{\lambda}+\dfrac{a(x/|x|)}{|x|^2}u_{\lambda}
+\lambda^2f(\lambda x, \lambda^2 t, u_{\lambda})=0\quad\text{in }\R^N\times (0,T/\lambda^2),
\end{equation}
in the sense that
\begin{align*}
&\int_\tau^{\frac{T}{\lambda^2}}\|u_{\lambda}(\cdot,t)\|^2_{{\mathcal
      H}_t}\,dt<+\infty,\quad\int_\tau^{\frac{T}{\lambda^2}}
\Big\|(u_{\lambda})_t+\frac{\nabla u_{\lambda}\cdot
    x}{2t}\Big\|^2_{({\mathcal H}_t)^\star}\!\!<+\infty \text{ for all
  }\tau\in \Big(0,{\frac{T}{\lambda^2}}\Big),\\
  &{\phantom{\bigg\langle}}_{{\mathcal
      H}_t^\star}\bigg\langle (u_{\lambda})_t+\frac{\nabla
    u_{\lambda}\cdot x}{2t},w
  \bigg\rangle_{{\mathcal H}_t}\\
  &\notag\qquad= \int_{\R^N}\bigg(\nabla u_{\lambda}(x,t)\cdot \nabla
  w(x)-
  \dfrac{a(x/|x|)}{|x|^2}\,u_{\lambda}(x,t)w(x)-\lambda^2 f(\lambda
  x, \lambda^2 t, u_{\lambda}(x,t))w(x)\bigg)G(x,t)\,dx
\end{align*}
for a.e. $t\in \big(0,{\frac{T}{\lambda^2}}\big)$ and for each $w\in
{\mathcal H}_t$. The frequency function associated to the scaled
equation (\ref{lambda}) is
\begin{equation}\label{eq:34}
{\mathcal N}_\lambda(t)=\frac{t\,D_\lambda(t)}{H_\lambda(t)},
\end{equation}
where
\begin{align*}
  D_{\lambda}(t)&=\int_{\R^N}\bigg(|\nabla u_{\lambda}(x,t)|^2-
  \frac{a(x/|x|)}{|x|^2}u_{\lambda}^2(x,t)-\lambda^2
  f(\lambda x, \lambda^2 t, u_{\lambda}(x,t))u_{\lambda}(x,t)\bigg)G(x, t)\,dx,\\
  H_{\lambda}(t)&=\int_{\R^N}u_{\lambda}^2(x,t) G(x,t)\,dx.
\end{align*}
The scaling properties of the operator combined with a suitable change
of variables easily imply that
\begin{align}\label{eq:19}
D_\lambda(t)=\lambda^2D(\lambda^2t)\quad\text{and}
\quad H_\lambda(t)=H(\lambda^2t),
\end{align}
and consequently
\begin{align}\label{eq:scale_for_N}
{\mathcal N}_\lambda(t)={\mathcal N}(\lambda^2t)\quad
\text{for all }t\in\Big(0,{\frac{2\alpha}{\lambda^2}}\Big).
\end{align}

\begin{Lemma}\label{l:blow_up}
  Let $a\in L^{\infty}\big({\mathbb S}^{N-1}\big)$ satisfy
  (\ref{eq:posde}) and $u\not\equiv 0$ be, in the sense of Definition
  \ref{def:solution}, either a weak solution to (\ref{prob1}), with
  $h$ satisfying (\ref{eq:der}) and (\ref{eq:h}), or a weak solution
  to (\ref{prob2}) satisfying (\ref{eq:u1}--\ref{eq:u2}) with
  $\varphi\in C^1(\R^N\times(0,T)\times\R)$ under assumption
  (\ref{eq:fi}). Let $\gamma:=\lim_{t\to 0^+}{\mathcal N}(t)$ as in
  Lemma \ref{l:limit}.  Then
\begin{itemize}
\item[(i)] $\gamma$ is an eigenvalue of the operator $L$ defined in
  (\ref{eq:13});
\item[(ii)] for every  sequence $\lambda_n\to 0^+$,
there exists a subsequence $\{\lambda_{n_k}\}_{k\in\N}$
and an eigenfunction $g$ of the operator $L$ associated to $\gamma$
such that, for all $\tau\in (0,1)$,
\begin{equation*}
\lim_{k\to+\infty}\int_\tau^1
\bigg\|\frac{u(\lambda_{n_k}x,\lambda_{n_k}^2t)}{\sqrt{H(\lambda_{n_k}^2)}}
-t^{\gamma}g(x/\sqrt t)\bigg\|_{{\mathcal H}_t}^2dt=0
\end{equation*}
and
$$
\lim_{k\to+\infty}\sup_{t\in[\tau,1]}
\bigg\|\frac{u(\lambda_{n_k}x,\lambda_{n_k}^2t)}{\sqrt{H(\lambda_{n_k}^2)}}
-t^{\gamma}g(x/\sqrt t)\bigg\|_{{\mathcal L}_t}=0.
$$

\end{itemize}
\end{Lemma}
\begin{pf}
Let
\begin{align}\label{eq:33}
w_{\lambda}(x,t):=\frac{u_\lambda(x,t)}{\sqrt{H(\lambda^2)}},
\end{align}
with $\lambda\in(0,\sqrt T_0)$, so that $1<T_0/\lambda^2$.  From Lemma
\ref{l:Hcreas} we obtain that, for all $t\in(0,1)$,
\begin{equation}\label{eq:21}
  \int_{\R^N}w_{\lambda}^2(x,t)G(x,t)\,dx
  =\frac{H(\lambda^2t)}{H(\lambda^2)}\leq t^{2C_1-\frac{N-2}{2}},
\end{equation}
with $C_1$ as in (\ref{eq:5}).
Lemma \ref{l:Nabove}, Corollaries \ref{c:pos_per} and
\ref{c:pos_per_nonlin}, and \eqref{eq:19} imply that
\begin{multline*}
  \frac1t \bigg(-\frac{N-2}4+{C_4}\bigg({\mathcal
    N}(T_0)+\frac{N-2}4\bigg)\bigg) H_\lambda(t)\geq
  \lambda^2D(\lambda^2 t)\\
  \geq \frac1t\bigg(C_1-\frac{N-2}{4}\bigg)H_\lambda(t)+C_1
  \int_{\R^N}|\nabla u_{\lambda}(x,t)|^2G(x,t)\,dx
\end{multline*}
and hence, in view of \eqref{eq:21},
\begin{align}\label{eq:20}
  t\int_{\R^N}|\nabla w_{\lambda}(x,t)|^2G(x,t)\,dx&\leq
 C_1^{-1} \big({\textstyle{C_4\big({\mathcal
    N}(T_0)+\frac{N-2}4\big)-C_1}}\big)
  \int_{\R^N}w_{\lambda}^2(x,t)G(x,t)\,dx\\
  &\notag \leq
C_1^{-1} \big({\textstyle{C_4\big({\mathcal
    N}(T_0)+\frac{N-2}4\big)-C_1}}\big)t^{2C_1-\frac{N-2}{2}},
\end{align}
for a.e. $t\in (0,1)$.
Let us consider the family of functions
$$
\widetilde{w}_{\lambda}(x,t)=w_{\lambda}(\sqrt{t}x,t)=\dfrac{u (\lambda \sqrt{t}x,
  \lambda^2t)}{\sqrt{H(\lambda^2)}},
$$
which, by scaling, satisfy
\begin{equation}\label{eq:22}
  \int_{\R^N} \widetilde{w}_{\lambda}^{2}(x,t)G(x,1)\,dx
  =\int_{\R^N} w_{\lambda}^{2}(x,t)G(x,t)\,dx
\end{equation}
and
\begin{equation}\label{eq:23}
  \int_{\R^N} |\nabla\widetilde{w}_{\lambda}(x,t)|^{2}G(x,1)\,dx
  =t\int_{\R^N}|\nabla  w_{\lambda}(x,t)|^{2}G(x,t)\,dx.
\end{equation}
From \eqref{eq:21}, \eqref{eq:20}, \eqref{eq:22}, and \eqref{eq:23},
we deduce that, for all  $\tau \in (0,1)$,
\begin{equation}\label{eq:27}
  \big\{\widetilde{w}_{\lambda}\big\}_{\lambda\in(0,\sqrt T_0)}
  \text{ is bounded in }
  L^\infty(\tau,1;{\mathcal H})
\end{equation}
uniformly with respect to $\lambda\in(0,\sqrt T_0)$.
Since
$$
\widetilde{w}_{\lambda}(x,t)=\frac{v(x,\lambda^2t)}{\sqrt{H(\lambda^2)}}
\quad\text{and}\quad
(\widetilde{w}_{\lambda})_t(x,t)=\frac{\lambda^2}{\sqrt{H(\lambda^2)}}\,
v_t(x,\lambda^2t)
$$
with $v$ as in Remark \ref{rem:uv}, from
(\ref{eq:24}) we deduce that, for all $\phi\in{\mathcal H}$,
\begin{multline}\label{eq:25}
  {\phantom{\big\langle}}_{{\mathcal H}^\star}\big\langle
  (\widetilde{w}_{\lambda})_t,\phi
  \big\rangle_{{\mathcal H}}
  =\frac1t\int_{\R^N}\bigg(\nabla \widetilde{w}_{\lambda}(x,t) \cdot
  \nabla \phi(x)-
  \dfrac{a\big(\frac{x}{|x|}\big)}{|x|^2}\,\widetilde{w}_{\lambda}(x,t)
  \phi(x)\\
  -\frac{\lambda^2t}{\sqrt{H(\lambda^2)}}f\Big(\lambda\sqrt tx, \lambda^2 t,
  \sqrt{H(\lambda^2)}\widetilde{w}_{\lambda}(x,t)\Big)\phi(x)\bigg)G(x,1)\,dx.
\end{multline}
In case {\bf (I)}, from (\ref{eq:h}) and Lemma \ref{Hardytemp}, we can estimate
the last term in the above integral as
\begin{align}\label{eq:26}
  &\lambda^2\left|\int_{\R^N}h(\lambda\sqrt
    tx, \lambda^2 t)\widetilde{w}_{\lambda}(x,t)\phi(x)G(x,1)\,dx\right|\\
  \notag&\leq
  C_h\lambda^2\int_{\R^N}|\widetilde{w}_{\lambda}(x,t)||\phi(x)|G(x,1)\,dx
  +C_h\frac{\lambda^\e}{t}\int_{\R^N}|x|^{-2+\e}
  |\widetilde{w}_{\lambda}(x,t)||\phi(x)|G(x,1)\,dx\\
  \notag& \leq
  C_h\lambda^2\|\widetilde{w}_{\lambda}(\cdot,t)\|_{\mathcal H}
  \|\phi\|_{\mathcal H}+ C_h\frac{\lambda^\e}{t}\int_{|x|\leq1} \frac{
    |\widetilde{w}_{\lambda}(x,t)||\phi(x)|}{|x|^2}G(x,1)\,dx\\
  \notag& \hskip6cm +C_h\frac{\lambda^\e}{t}\int_{|x|\geq1}
  |\widetilde{w}_{\lambda}(x,t)||\phi(x)|G(x,1)\,dx\\
  \notag&\leq
  C_h\frac{\lambda^\e}{t}\bigg(t\,\lambda^{2-\e}+\frac{\max\{4,N-2\}}{(N-2)^2}+1
  \bigg) \|\widetilde{w}_{\lambda}(\cdot,t)\|_{\mathcal H}
  \|\phi\|_{\mathcal H}
\end{align}
for all $\lambda\in (0,\sqrt T_0)$ and a.e. $t\in(0,1)$. From
(\ref{eq:25}), (\ref{eq:26}), and Lemma \ref{Hardytemp} it follows that,
for all $\lambda\in (0,\sqrt T_0)$ and a.e. $t\in(0,1)$,
\begin{multline*}
  \big| \!\!{\phantom{\big\langle}}_{{\mathcal H}^\star}\big\langle
  (\widetilde{w}_{\lambda})_t,\phi \big\rangle_{{\mathcal
      H}}\big|\\\leq
  \Big({\textstyle{1+\frac{\max\{4,N-2\}}{(N-2)^2}\|a\|_{L^{\infty}({\mathbb
          S}^{N-1})}
      +C_hT_0^{\e/2}\Big(T_0^{1-\e/2}+\frac{\max\{4,N-2\}}{(N-2)^2}+1
      \Big)}}\Big)\frac{\|\widetilde{w}_{\lambda}(\cdot,t)\|_{\mathcal
      H} \|\phi\|_{\mathcal H}}{t}
\end{multline*}
and hence
\begin{equation}\label{eq:89}
\|(\widetilde{w}_{\lambda})_t(\cdot,t)\|_{{\mathcal H}^\star}\leq
\frac{{\rm const}}{t}\|\widetilde{w}_{\lambda}(\cdot,t)\|_{\mathcal H}.
\end{equation}
In case {\bf (II)}, from (\ref{eq:fi}), H\"older's inequality, and Lemma
\ref{l:sob}, we obtain
\begin{multline}\label{eq:90}
  \bigg|\frac{\lambda^2}{\sqrt{H(\lambda^2)}}\int_{\R^N}\varphi\Big(\lambda\sqrt
  tx, \lambda^2 t,
  \sqrt{H(\lambda^2)}\widetilde{w}_{\lambda}(x,t)\Big)\phi(x)G(x,1)\,dx
  \bigg|\\
  \leq C_\varphi\frac{\lambda^2}{\sqrt{H(\lambda^2)}}
  \int_{\R^N}\Big(\sqrt{H(\lambda^2)}|\widetilde{w}_{\lambda}(x,t)|+
  (\sqrt{H(\lambda^2)})^p|\widetilde{w}_{\lambda}(x,t)|^p\Big)
  |\phi(x)|G(x,1)\,dx\\
  \leq C_\varphi
  \lambda^2\int_{\R^N}|\widetilde{w}_{\lambda}(x,t)||\phi(x)|G(x,1)\,dx+
  C_\varphi
  \lambda^2(H(\lambda^2))^{\frac{p-1}{2}}\int_{\R^N}|\widetilde{w}_{\lambda}(x,t)|^p
  |\phi(x)|G(x,1)\,dx\\
  \leq \|\widetilde{w}_{\lambda}(\cdot,t)\|_{\mathcal H}
  \|\phi\|_{\mathcal H}
\frac{\lambda^{\frac{N+2-p(N-2)}{p+1}}}{t}
\bigg(C_\varphi t
  \lambda^{\frac{N(p-1)}{p+1}}+ C_\varphi C_{p+1}
  t^{\frac{N+2-p(N-2)}{2(p+1)}}
  \bigg(\int_{\R^N}|u(x,\lambda^2t)|^{p+1}\,dx\bigg)^{\!\!\frac{p-1}{p+1}}\bigg).
\end{multline}
From (\ref{eq:25}), (\ref{eq:90}), Lemma \ref{Hardytemp}, the fact
that $p<2^*-1$, and (\ref{eq:u1}), it follows that, for all
$\lambda\in (0,\sqrt T_0)$ and a.e. $t\in(0,1)$, estimate
(\ref{eq:89}) holds also in case {\bf (II)}.
Then, in view of (\ref{eq:27}),  estimate (\ref{eq:89}) yields, for all
 $\tau \in (0,1)$,
\begin{equation}\label{eq:28}
  \big\{(\widetilde{w}_{\lambda})_t\big\}_{\lambda\in(0,\sqrt T_0)}
  \text{ is bounded in }
  L^\infty(\tau,1;{\mathcal H}^\star)
\end{equation}
uniformly with respect to $\lambda\in(0,\sqrt T_0)$.  From
(\ref{eq:27}), (\ref{eq:28}), and \cite[Corollary 8]{S}, we deduce
that $\big\{\widetilde{w}_{\lambda}\big\}_{\lambda\in(0,\sqrt T_0)}$  is relatively
compact in $C^0([\tau,1],{\mathcal L})$ for all $\tau\in (0,1)$.
 Therefore, for any given sequence $\lambda_n\to 0^+$,
there exists a subsequence $\lambda_{n_k}\to0^+$ such that
\begin{equation}\label{eq:29}
  \widetilde{w}_{\lambda_{n_k}}\to\widetilde{w} \quad\text{in} \quad
  C^0([\tau,1],{\mathcal L})
\end{equation}
for all $\tau\in(0,1)$ and for some $\widetilde{w}\in
\bigcap_{\tau\in(0,1)}C^0([\tau,1],{\mathcal L})$.  We notice that a
diagonal procedure allows subtracting a subsequence which does not
depend on $\tau$.  Since
$$
1=\|\widetilde{w}_{\lambda_{n_k}}(\cdot,1)\|_{\mathcal L},
$$
 the convergence (\ref{eq:29}) ensures that
\begin{equation}\label{eq:41}
\|\widetilde{w}(\cdot,1)\|_{\mathcal L}=1.
\end{equation}
In particular $\widetilde{w}$ is nontrivial. Furthermore, by
(\ref{eq:27}) and (\ref{eq:28}), the subsequence can be chosen in such a way
that also
\begin{equation}\label{eq:30}
  \widetilde{w}_{\lambda_{n_k}}\weakly\widetilde{w} \quad\text{weakly in }
  L^2(\tau,1;{\mathcal H})
  \quad\text{and}\quad
  (\widetilde{w}_{\lambda_{n_k}})_t\weakly\widetilde{w}_t
  \quad\text{weakly in }
  L^2(\tau,1;{\mathcal H}^\star)
\end{equation}
for all $\tau\in(0,1)$; in particular $\widetilde{w}\in
\bigcap_{\tau\in(0,1)}L^2(\tau,1;{\mathcal H})$ and
$\widetilde{w}_t\in \bigcap_{\tau\in(0,1)}L^2(\tau,1;{\mathcal
  H}^\star)$. We now claim that
\begin{equation}\label{eq:31}
  \widetilde{w}_{\lambda_{n_k}}\to\widetilde{w} \quad\text{strongly in} \quad
  L^2(\tau,1;{\mathcal H})\quad\text{for all }\tau\in(0,1).
\end{equation}
To prove the claim, we notice that  (\ref{eq:30}) allows
passing to the limit in (\ref{eq:25}). Therefore, in view of
(\ref{eq:26}) and (\ref{eq:89}) which ensure the vanishing at the
limit of the perturbation term,
\begin{equation}\label{eq:32}
  {\phantom{\big\langle}}_{{\mathcal
      H}^\star}\big\langle \widetilde{w}_t,\phi
  \big\rangle_{{\mathcal H}}=\frac1t\int_{\R^N}\bigg(\nabla \widetilde{w}(x,t)
  \cdot \nabla \phi(x)-
  \dfrac{a\big(\frac{x}{|x|}\big)}{|x|^2}\,\widetilde{w}(x,t)
  \phi(x)\bigg)G(x,1)\,dx
\end{equation}
for all $\phi\in{\mathcal H}$ and a.e. $t \in (0,1)$, i.e. $\widetilde
w$ is a weak solution to
\begin{equation*}
  \widetilde{w}_t+\frac1t\bigg(\Delta \widetilde{w}-
  \frac x2\cdot \nabla \widetilde{w}+
  \dfrac{a(x/|x|)}{|x|^2}\,\widetilde{w}\bigg)=0.
\end{equation*}
Testing the difference between (\ref{eq:25}) and (\ref{eq:32}) with
$(\widetilde{w}_{\lambda_{n_k}}- \widetilde{w})$ and integrating
with respect to $t$ between $\tau$ and $1$, we obtain
\begin{multline*}
  \int_\tau^1\bigg(\int_{\R^N}\bigg(|\nabla (\widetilde{w}_{\lambda_{n_k}}-
  \widetilde{w})(x,t)|^2-\frac{a(x/|x|)}{|x|^2}\,
  |(\widetilde{w}_{\lambda_{n_k}}- \widetilde{w})(x,t)|^2\bigg)
  \,G(x,1)\,dx\bigg)dt\\
  = \frac{1}2\|\widetilde{w}_{\lambda_{n_k}}(1)-
  \widetilde{w}(1)\|^2_{\mathcal L}-
  \frac{\tau}2\|\widetilde{w}_{\lambda_{n_k}}(\tau)-
  \widetilde{w}(\tau)\|^2_{\mathcal L}-\int_\tau^1\bigg(\int_{\R^N}
  |(\widetilde{w}_{\lambda_{n_k}}-
  \widetilde{w})(x,t)|^2  \,G(x,1)\,dx\bigg)dt\\
  +\frac{\lambda_{n_k}^2}{\sqrt{H(\lambda_{n_k}^2)}}\int_\tau^1\bigg(\int_{\R^N}
  t f\Big(\lambda_{n_k}\sqrt tx, \lambda_{n_k}^2 t,
  \sqrt{H(\lambda_{n_k}^2)}\widetilde{w}_{\lambda_{n_k}}(x,t)\Big)
  (\widetilde{w}_{\lambda_{n_k}}- \widetilde{w})(x,t) G(x,1)\,dx \bigg)dt.
\end{multline*}
Then, from (\ref{eq:26}), (\ref{eq:90}), and (\ref{eq:29}), we obtain that,
for all $\tau\in (0,1)$,
$$
\lim_{k\to+\infty} \int_\tau^1\bigg(\int_{\R^N}\bigg(|\nabla
(\widetilde{w}_{\lambda_{n_k}}-
\widetilde{w})(x,t)|^2-\frac{a(x/|x|)}{|x|^2}\,
|(\widetilde{w}_{\lambda_{n_k}}- \widetilde{w})(x,t)|^2\bigg)
\,G(x,1)\,dx\bigg)dt=0,
$$
which, by Corollary \ref{c:pos_def} and (\ref{eq:29}), implies
the convergence claimed in (\ref{eq:31}).
Thus, we have obtained that, for all $\tau\in(0,1)$,
\begin{equation}\label{eq:35}
\lim_{k\to+\infty}\int_\tau^1
\|{w}_{\lambda_{n_k}}(\cdot,t)-w(\cdot,t)\|_{{\mathcal H}_t}^2dt=0
\end{equation}
and
$$
\lim_{k\to+\infty}\sup_{t\in[\tau,1]}
\|{w}_{\lambda_{n_k}}(\cdot,t)-w(\cdot,t)\|_{{\mathcal L}_t}=0,
$$
where
$$
w(x,t):=\widetilde{w}\Big(\frac{x}{\sqrt t },t\Big)
$$
is a weak solution (in the sense of Definition \ref{def:solution})
of
\begin{equation}\label{eq:limit_equation}
w_t+\Delta w+\dfrac{a(x/|x|)}{|x|^2}\,w=0.
\end{equation}
We notice that, by (\ref{eq:34}) and (\ref{eq:33}),
\begin{align*}
  &{\mathcal N}_\lambda(t)\\
  &=\frac{ t\int_{\R^N}\Big(|\nabla w_{\lambda}(x,t)|^2-
    \frac{a(x/|x|)}{|x|^2}w_{\lambda}^2(x,t)-\frac{\lambda^2}{\sqrt{H(\lambda^2)}}
    f\big(\lambda x,\lambda^2 t, \sqrt{H(\lambda^2)}w_{\lambda}(x,t)\big)
    w_{\lambda}(x,t)\Big)G(x, t)\,dx}{\int_{\R^N}w_{\lambda}^2(x,t)
    G(x,t)\,dx}
\end{align*}
for all $t\in(0,1)$.  Since, by \eqref{eq:35},
${w}_{\lambda_{n_k}}(\cdot,t)\to w(\cdot,t)$ in ${\mathcal H}_t$ for
a.e. $t\in (0,1)$, and, by \eqref{eq:26} and~(\ref{eq:90}),
$$
\frac{t\lambda^2_{n_k}}{\sqrt{H(\lambda^2_{n_k})}}\int_{\R^N}
    f\Big(\lambda_{n_k} x,\lambda_{n_k}^2 t, \sqrt{H(\lambda_{n_k}^2)}
w_{\lambda_{n_k}}(x,t)\Big)
    w_{\lambda_{n_k}}(x,t)G(x, t)\,dx
\to 0
$$
for a.e. $t\in (0,1)$, we obtain that
\begin{multline}\label{eq:43}
  \int_{\R^N}\bigg(|\nabla w_{\lambda_{n_k}}(x,t)|^2-
  \frac{a(x/|x|)}{|x|^2}w_{\lambda_{n_k}}^2(x,t)-\\
  \frac{\lambda^2_{n_k}}{\sqrt{H(\lambda^2_{n_k})}}
  f\Big(\lambda_{n_k} x,\lambda^2 t,
  \sqrt{H(\lambda^2_{n_k})}w_{\lambda_{n_k}}(x,t)\Big)
  w_{\lambda_{n_k}}(x,t)\bigg)G(x, t)\,dx \to D_w(t)
\end{multline}
and
\begin{equation}\label{eq:44}
\int_{\R^N}w_{\lambda_{n_k}}^2(x,t)
  G(x,t)\,dx\to H_w(t)
\end{equation}
for a.e. $t\in (0,1)$, where
$$
D_w(t)=\int_{\R^N}\bigg(|\nabla w(x,t)|^2-
  \frac{a(x/|x|)}{|x|^2}w^2(x,t)\bigg)G(x, t)\,dx
\ \text{ and }\
H_w(t)=\int_{\R^N}w^2(x,t)
  G(x,t)\,dx.
$$
We point out that
\begin{equation}\label{eq:42}
H_w(t)>0 \quad\text{for all}\quad t\in(0,1);
\end{equation}
indeed, (\ref{eq:41})
yields
\begin{equation}\label{eq:wl1}
  \int_{\R^N}w^2(x,1)G(x,1)\,dx=1,
\end{equation}
 which,
arguing as in Lemma \ref{l:t_0=0} or applying directly the
Unique Continuation Principle proved by \cite[Theorem 1.2]{poon} to equation
(\ref{eq:limit_equation}),  implies
that $\int_{\R^N}w^2(x,t)
  G(x,t)\,dx>0$ for all $t\in(0,1)$.
From (\ref{eq:43}) and (\ref{eq:44}), it follows that
\begin{equation}\label{eq:36}
{\mathcal N}_{\lambda_{n_k}}(t)\to {\mathcal N}_w(t)
\quad\text{for a.e. }t\in(0,1),
\end{equation}
where ${\mathcal N}_w$ is the frequency function associated to the limit
equation \eqref{eq:limit_equation}, i.e.
\begin{equation}\label{eq:Nw}
{\mathcal N}_w(t)=\frac{ tD_w(t)}
{H_w(t)},
\end{equation}
which is well defined on $(0,1)$ by (\ref{eq:42}).

On the other hand, \eqref{eq:scale_for_N} implies that ${\mathcal
  N}_{\lambda_{n_k}}(t)={\mathcal N}(\lambda^2_{n_k}t)$ for all $t\in
(0,1)$ and $k\in\N$.  Fixing $t\in (0,1)$ and passing to the limit as
$k\to+\infty$, from Lemma \ref{l:limit} we obtain
\begin{equation}\label{eq:37}
{\mathcal N}_{\lambda_{n_k}}(t)\to \gamma
\quad\text{for all }t\in(0,1).
\end{equation}
Combining \eqref{eq:36} and \eqref{eq:37}, we deduce that
\begin{equation}\label{eq:Nw-gamma}
{\mathcal N}_w(t)=\gamma\quad \text{for all }t\in(0,1).
\end{equation}
Therefore ${\mathcal N}_w$ is constant in $(0,1)$ and hence ${\mathcal
  N}_w'(t)=0$ for any $t\in (0,1)$.  By \eqref{eq:limit_equation} and
Lemma~\ref{l:Nprime} with $f\equiv 0$, we obtain
\begin{multline*}
\bigg(\int_{\R^N}\Big|w_t(x,t)+\frac{\nabla w(x,t)\cdot
        x}{2t}\Big|^2G(x,t)\,dx\bigg)
      \bigg(\int_{\R^N}w^2(x,t)\, G(x,t)\,dx\bigg)\\
-\bigg(\int_{\R^N}\Big(w_t(x,t)+\frac{\nabla
        w(x,t)\cdot x}{2t}\Big)w(x,t)G(x,t)\,dx
      \bigg)^{\!2}=0 \quad
  \text{for all } t\in (0,1),
\end{multline*}
i.e.
$$
\left(w_t(\cdot,t)+\frac{\nabla w(\cdot,t)\cdot
    x}{2t},w(\cdot,t)\right)^2_{{\mathcal L}_t}
=\left\|w_t(\cdot,t)+\frac{\nabla w(\cdot,t)\cdot
    x}{2t}\right\|^2_{{\mathcal L}_t} \|w(\cdot,t)\|^2_{{\mathcal
    L}_t},
$$
where $(\cdot,\cdot)_{{\mathcal L}_t}$ denotes the scalar product in
${\mathcal L}_t$.  This shows that, for all $t\in(0,1)$,
$w_t(\cdot,t)+\frac{\nabla w(\cdot,t)\cdot x}{2t}$ and $w(\cdot,t)$
have the same direction as vectors in ${\mathcal L}_t$ and hence there
exists a  function $\beta:(0,1)\to\R$ such that
\begin{equation}\label{eq:38}
w_t(x,t)+\frac{\nabla w(x,t)\cdot x}{2t} =\beta(t)w(x,t)\quad\text{for
  a.e. }t\in(0,1)\text{ and a.e. }x\in\R^N.
\end{equation}
Testing \eqref{eq:limit_equation} with $\phi=w(\cdot,t)$ in the sense
of \eqref{eq:defsol2} and taking into account \eqref{eq:38},
we find that
\begin{align*}
D_w(t)=
{\phantom{\bigg\langle}}_{{\mathcal
    H}_t^\star}\bigg\langle
w_t(\cdot,t)+\frac{\nabla w(\cdot,t)\cdot x}{2t},w(\cdot,t)
\bigg\rangle_{{\mathcal H}_t}=\beta(t)H_w(t),
\end{align*}
which, by (\ref{eq:Nw}) and (\ref{eq:Nw-gamma}), implies that
$$
\beta(t)=\frac{\gamma}{t}\quad\text{for a.e. }t\in(0,1).
$$
Hence (\ref{eq:38}) becomes
\begin{equation}\label{eq:45}
w_t(x,t)+\frac{\nabla w(x,t)\cdot x}{2t}
=\frac{\gamma}{t}\,w(x,t)\text{ for a.e. }(x,t)\in\R^N\times
(0,1)\text{ and in a distributional sense}.
\end{equation}
Combining (\ref{eq:45}) with (\ref{eq:limit_equation}), we obtain
\begin{equation}\label{eq:46}
\Delta w+\dfrac{a(x/|x|)}{|x|^2}\,w
-\frac{\nabla w(x,t)\cdot x}{2t}+
\frac{\gamma}{t}\,w(x,t)=0
\end{equation}
for a.e. $(x,t)\in\R^N\times
(0,1)$  and in a weak  sense.
From (\ref{eq:45}), it follows that, letting, for all $\eta>0$ and a.e.
$(x,t)\in\R^N\times(0,1)$, $w^\eta(x,t):=w(\eta x,\eta^2t)$,
there holds
$$
\frac{d
  w^\eta}{d\eta}=\frac{2\gamma}{\eta}w^\eta
$$
a.e. and in a
distributional sense. By integration, we obtain that
\begin{equation}
  \label{eq:weta}
 w^\eta(x,t)=w(\eta x,\eta^2t)=\eta^{2\gamma}w(x,t) \quad\text{for all }\eta>0
\text{ and a.e. }(x,t)\in\R^N\times(0,1).
\end{equation}
Let
$$
g(x)=w(x,1);
$$
from (\ref{eq:wl1}), we have that $g\in{\mathcal L}$, $\|g\|_{{\mathcal L}}=1$,
 and,
from (\ref{eq:weta}),
\begin{equation}\label{eq:47}
w(x,t)=w^{\sqrt t }\Big(\frac{x}{\sqrt t
},1\Big)=t^{\gamma}w\Big(\frac{x}{\sqrt t },1\Big)=
t^{\gamma}g\Big(\frac{x}{\sqrt t}\Big)
\quad\text{for a.e. }(x,t)\in\R^N\times(0,1).
\end{equation}
In particular, from (\ref{eq:47}), $g\big({\cdot}/{\sqrt t}\big)\in
{\mathcal H}_t$ for a.e. $t\in (0,1)$ and hence, by scaling,
$g\in{\mathcal H}$.
From (\ref{eq:46}) and (\ref{eq:47}), we obtain
that $g\in{\mathcal H}\setminus\{0\}$ weakly solves
$$
-\Delta g(x)+
\frac{\nabla g(x)\cdot x}{2}-\frac{a(x/|x|)}{|x|^2}\,g(x)=\gamma\, g(x),
$$
i.e.
$\gamma$ is an eigenvalue of the operator $L$ defined in (\ref{eq:13})
and $g$ is an eigenfunction of $L$ associated to $\gamma$. The proof
is now complete.
\end{pf}

\noindent Let us now describe the behavior of $H(t)$ as $t\to 0^+$.
\begin{Lemma} \label{l:limite} Under the same assumptions as in Lemma
  \ref{l:blow_up}, let $\gamma:=\lim_{t\rightarrow 0^+} {\mathcal N}(t)$ be as in
  Lemma \ref{l:limit}.  Then the limit
\[
\lim_{t\to 0^+}t^{-2\gamma}H(t)
\]
exists and it is finite.
\end{Lemma}
\begin{pf}
In view of \eqref{eq:52}, it is sufficient to prove that the limit exists.
By \eqref{eq:10i},  Lemma \ref{l:limit}, and Lemma \ref{l:Nprime},  we have,
for all $t\in(0,T_0)$,
\begin{align*}
  \frac{d}{dt} \frac{H(t)}{t^{2\gamma}} &=-2\gamma t^{-2\gamma-1}
  H(t)+t^{-2\gamma} H'(t) =2t^{-2\gamma-1} (tD(t)-\gamma H(t))\\
&=
  2t^{-2\gamma-1} H(t) \int_0^t (\nu_{1}(s)+\nu_{2}(s))\, ds,
\end{align*}
with $\nu_1,\nu_2$ as in (\ref{eq:nu1nu2}).
After integration over $(t,T_0)$,
\begin{equation}\label{inte}
  \frac{H(T_0)}{T_0^{2\gamma}}-\frac{H(t)}{t^{2\gamma}}=
\int_t^{T_0} 2s^{-2\gamma-1}
  H(s) \left( \int_0^s \nu_{1}(r)dr \right) ds +
\int_t^{T_0} 2s^{-2\gamma-1}
  H(s) \left( \int_0^s \nu_{2}(r)dr \right) ds.
\end{equation}
By \eqref{eq:17},  $\nu_{1}(t)\geq 0$ and hence
$$
\lim_{t\to 0^+} \int_t^{T_0} 2s^{-2\gamma-1} H(s) \left( \int_0^s
  \nu_{1}(r)dr \right) ds
$$
exists.  On the other hand, by Lemmas \ref{l:est_N_2} and
\ref{l:Nabove} we have that $s^{-\delta}\int_0^s|\nu_{2}(r)|dr$ is bounded in
$(0,T_0)$ with $\delta$ defined in (\ref{eq:delta}), while, from
Lemma \ref{stimaH}, we deduce that $t^{-2\gamma}H(t)$  is bounded in
$(0,T_0)$. Therefore, for some ${\rm const}>0$, there holds
$$
\left|
2s^{-2\gamma-1}
  H(s) \left( \int_0^s \nu_{2}(r)dr \right)
\right|\leq {\rm const\,} s^{-1+\delta}
$$
for all $s\in(0,T_0)$, which proves that $s^{-2\gamma-1} H(s)
\left( \int_0^s \nu_{2}(r) dr \right)\in L^1(0,T_0)$.  We
conclude that both terms at the right hand side of (\ref{inte}) admit
a limit as $t\to 0^+$ thus completing the proof.
\end{pf}

\noindent In the following lemma, we prove that  $\lim_{t\to
  0^+}t^{-2\gamma}H(t)$ is indeed strictly positive.

\begin{Lemma}\label{limite_pos}
  Under the same assumptions as in Lemma \ref{l:blow_up} and letting
  $\gamma:=\lim_{t\rightarrow 0^+} {\mathcal N}(t)$ be as in Lemma
  \ref{l:limit}, there holds
\[
\lim_{t\to 0^+}t^{-2\gamma}H(t)>0.
\]
\end{Lemma}
\begin{pf}
  Let us assume by contradiction that $\lim_{t\to
    0^+}t^{-2\gamma}H(t)=0$ and let $\{\widetilde V_{n,j}:
  j,n\in\N,j\geq 1\}$ be the orthonormal basis of ${\mathcal L}$
  introduced in Remark \ref{rem:ortho}.  Since
  $u_{\lambda}(x,1)=u(\lambda x,\lambda^2)\in{\mathcal L}$ for all
  $\lambda\in (0,\sqrt T_0)$, $u_{\lambda}(x,1)\in{\mathcal H}$ for
  a.e.  $\lambda\in (0,\sqrt T_0)$, and $f(\lambda x,\lambda^2,
  u_{\lambda}(x,1))\in{\mathcal H}^\star$ for a.e. $\lambda\in
  (0,\sqrt T_0)$, we can expand them as
\begin{align}\label{eq:70}
  u_{\lambda}(x,1)&=\sum_{\substack{m,k\in\N\\k\geq1}}
  u_{m,k}(\lambda)
  \widetilde V_{m,k}(x)\quad \text{in }{\mathcal L},\\
  \notag f(\lambda x,\lambda^2,
  u_{\lambda}(x,1))&=\sum_{\substack{m,k\in\N\\k\geq1}}
  \xi_{m,k}(\lambda) \widetilde V_{m,k}(x)\quad \text{in }{\mathcal
    H}^\star,
\end{align}
where
\begin{equation}\label{eq:68}
u_{m,k}(\lambda)=\int_{\R^N}u_{\lambda}(x,1)\widetilde V_{m,k}(x)G(x,1)\,dx
\end{equation}
and
\begin{equation}\label{eq:69}
  \xi_{m,k}(\lambda)=
  {\phantom{\bigg\langle}}_{{\mathcal
      H}^\star}\bigg\langle
  f(\lambda \cdot,\lambda^2,u_{\lambda}(\cdot,1)), \widetilde V_{m,k}
  \bigg\rangle_{{\mathcal H}}
  =\int_{\R^N}f(\lambda x,\lambda^2,
  u_{\lambda}(x,1))\widetilde V_{m,k}(x)G(x,1)\,dx.
\end{equation}
By orthogonality of the $\widetilde V_{m,k}$'s in ${\mathcal L}$,
we have that
$$
H(\lambda^2)=\sum_{\substack{n,j\in\N\\j\geq1}} (u_{n,j}(\lambda) )^2
\geq (u_{m,k}(\lambda) )^2\quad\text{for all }
\lambda\in(0,\sqrt{T_0})\text{ and }m,k\in\N,\ k\geq1.
$$
Hence,  $\lim_{t\to
    0^+}t^{-2\gamma}H(t)=0$ implies that
\begin{equation}\label{eq:58}
\lim_{\lambda\to 0^+}\lambda^{-2\gamma}u_{m,k}(\lambda)=0
\quad\text{for all }
m,k\in\N,\ k\geq1.
\end{equation}
Moreover, we can show that the function $\lambda\mapsto u_{m,k}(\lambda)$
is absolutely continuous in $(0,\sqrt{T_0})$ and $u_{m,k}'(\lambda)=
{\phantom{\langle}}_{{\mathcal
    H}^\star}\langle
\frac{d}{d\lambda}u_{\lambda}(x,1),\widetilde V_{m,k}(x)
\rangle_{{\mathcal H}}$. Hence
$$
\frac{d}{d\lambda}u_{\lambda}(x,1)=
\sum_{\substack{m,k\in\N\\k\geq1}} u_{m,k}'(\lambda)
\widetilde V_{m,k}(x)\quad \text{in }{\mathcal H}^\star.
$$
Furthermore,
$$
\Delta u_{\lambda}(x,1)=\lambda^2\Delta u(\lambda x,\lambda^2)
=
\sum_{\substack{m,k\in\N\\k\geq1}} u_{m,k}(\lambda)
\Delta\widetilde V_{m,k}(x)\quad \text{in }{\mathcal H}^\star.
$$
From \eqref{prob} and the fact that $\widetilde V_{m,k}(x)$ is an
eigenfuntion of the operator $L$ associated to the eigenvalue
$\gamma_{m,k}$ defined in (\ref{eq:65}), it follows that
\begin{align*}
  \dfrac{d}{d\lambda}&u_{\lambda}(x,1)=2\lambda u_{t}(\lambda x,\lambda^2)+
  \nabla u(\lambda x,\lambda^2)\cdot x\\
  &=2\lambda \bigg(-\Delta u(\lambda x, \lambda^2)-
  \frac{a(x/|x|)}{\lambda^2|x|^2}u(\lambda x, \lambda^2)- f(\lambda
  x,\lambda^2,u(\lambda x, \lambda^2))\bigg)+\nabla u(\lambda
  x,\lambda^2)
  \cdot x\\
  &=\dfrac{2}{\lambda}\sum_{\substack{m,k\in\N\\k\geq1}}
  u_{m,k}(\lambda)\bigg(-\Delta \widetilde
  V_{m,k}(x)-\frac{a(x/|x|)}{|x|^2}\widetilde V_{m,k}(x) +\frac{\nabla
    \widetilde V_{m,k}\cdot x}2\bigg)
  -2\lambda\sum_{\substack{m,k\in\N\\k\geq1}}
  \xi_{m,k}(\lambda) \widetilde V_{m,k}(x)\\
  &=\dfrac{2}{\lambda}\sum_{\substack{m,k\in\N\\k\geq1}}
  \gamma_{m,k}u_{m,k}(\lambda)\widetilde V_{m,k}(x)-2\lambda
  \sum_{\substack{m,k\in\N\\k\geq1}} \xi_{m,k}(\lambda) \widetilde
  V_{m,k}(x).
\end{align*}
Therefore, we have that
$$
u'_{m,k}(\lambda)=\frac{2}{\lambda}\gamma_{m,k}u_{m,k}(\lambda)-2\lambda
\xi_{m,k}(\lambda) \quad\text{for all }m,k\in\N,\ k\geq 1,
$$
a.e. and distributionally in $(0,\sqrt{T_0})$.
By integration, we obtain, for all
$\lambda,\bar\lambda\in(0,\sqrt{T_0})$,
\begin{equation}\label{eq:59}
u_{m,k}(\bar\lambda)=
\bar\lambda^{2\gamma_{m,k}}\left(\lambda^{-2\gamma_{m,k}} u_{m,k}(\lambda)+
2\int_{\bar\lambda}^{\lambda}s^{1-2\gamma_{m,k}}\xi_{m,k}(s)\,ds\right).
\end{equation}
From Lemma \ref{l:blow_up}, $\gamma$ is an eigenvalue of the operator
$L$, hence, by Proposition \ref{p:explicit_spectrum}, there exist
$m_0,k_0\in\N$, $k_0\geq 1$, such that
$\gamma=\gamma_{m_0,k_0}=m_0-\frac{\alpha_{k_0}}2$.  Let us denote as
$E_0$ the associated eigenspace and by $J_0$ the finite set of indices
$\{(m,k)\in\N\times(\N\setminus\{0\}):\gamma=m-\frac{\alpha_{k}}2\}$,
so that $\#J_0=m(\gamma)$, with $m(\gamma)$ as in \eqref{eq:cardinal},
and an orthonormal basis of $E_0$ is given by $\{\widetilde V_{m,k} :
(m,k)\in J_0\}$. In order to estimate $\xi_{m,k}$, we distinguish between
 case {\bf (I)} and case {\bf (II)}.

\begin{description}
\item[Case {\bf (I)}] From \eqref{eq:h}, for all $(m,k)\in J_0$, we
  can estimate $\xi_{m,k}$ as
\begin{align}\label{eq:54}
&  |\xi_{m,k}(\lambda)|\leq
  C_h\int_{\R^N}(1+\lambda^{-2+\e}|x|^{-2+\e})|u(\lambda x,\lambda^2)|
  |\widetilde V_{m,k}(x)|G(x,1)\,dx\\
  \notag&\leq C_h\bigg(\int_{\R^N}u^2(\lambda x,\lambda^2)
  G(x,1)\,dx\bigg)^{\!\!1/2}\bigg(\int_{\R^N}\widetilde V_{m,k}^2(x)
  G(x,1)\,dx\bigg)^{\!\!1/2}\\
  \notag&\qquad+ C_h\lambda^{-2+\e/2}\int_{|x|\leq
    \lambda^{-1/2}}\frac{|u(\lambda x,\lambda^2)| |\widetilde
    V_{m,k}(x)|}{|x|^2}G(x,1)\,dx\\
  \notag&\qquad+C_h\lambda^{-1+\e/2}\int_{|x|\geq \lambda^{-1/2}}|u(\lambda
  x,\lambda^2)|
  |\widetilde V_{m,k}(x)|G(x,1)\,dx\\
  \notag&\hskip-0.5cm\leq C_h(1+\lambda^{-1+\frac\e2})\sqrt{H(\lambda^2)}
  +C_h\lambda^{-2+\frac\e2}\bigg(\int_{\R^N}{\textstyle{\frac{u^2(\lambda
    x,\lambda^2)} {|x|^2}}}G(x,1)\,dx\bigg)^{\!\!\frac12}
  \bigg(\int_{\R^N}{\textstyle{\frac{\widetilde V_{m,k}^2(x)}
  {|x|^2}G(x,1)\,dx}}\bigg)^{\!\!\frac12}.
\end{align}
From Corollary \ref{c:pos_per} and Lemma \ref{l:Nabove}, it follows that
\begin{multline}\label{eq:55}
  \int_{\R^N}\frac{u^2(\lambda x,\lambda^2)} {|x|^2}G(x,1)\,dx
  =\lambda^2\int_{\R^N}\frac{u^2(y,\lambda^2)} {|y|^2}G(y,\lambda^2)\,dy
  \leq \frac{\lambda^2}{C_1'}\bigg(D(\lambda^2)+
  \frac{C_2}{\lambda^2}H(\lambda^2)\bigg)\\=
  \frac{H(\lambda^2)}{C_1'}\big({\mathcal N}(\lambda^2)+C_2\big)\leq\frac{C_2
-\frac{N-2}{4}+C_{4}\big({\mathcal N}(T_0)+\frac{N-2}4\big)
}{C_1'}
H(\lambda^2),
\end{multline}
while, from Lemma \ref{Hardy_aniso}, for all $(m,k)\in J_0$,
\begin{equation}\label{eq:56}
\int_{\R^N}\frac{\widetilde V_{m,k}^2(x)}
  {|x|^2}G(x,1)\,dx\leq
\bigg(\mu_1(a)+\frac{(N-2)^2}4\bigg)^{-1}\bigg(\gamma+\frac{N-2}4
\bigg).
\end{equation}
From (\ref{eq:54}), (\ref{eq:55}), (\ref{eq:56}),  and Lemma \ref{stimaH},
we deduce that
\begin{align}\label{eq:57}
|\xi_{m,k}(\lambda)|\leq C_6 \lambda^{-2+\frac{\e}{2}+2\gamma},
\quad\text{for all }\lambda\in(0,\sqrt{T_0})
\end{align}
and for some positive constant $C_6$ depending on
$a,N,\gamma,h,T_0,K_1,\e$ but independent of $\lambda$ and $(m,k)\in
J_0$.

\medskip
\item[Case {\bf (II)}] From \eqref{eq:fi} and Lemma \ref{l:sob}, for
  all $(m,k)\in J_0$, we can estimate $\xi_{m,k}$ as
\begin{align}\label{eq:94}
  & |\xi_{m,k}(\lambda)|\leq C_\varphi\int_{\R^N}\big(|u(\lambda
  x,\lambda^2)|+|u(\lambda x,\lambda^2)|^p\big)
  |\widetilde V_{m,k}(x)|G(x,1)\,dx\\
  \notag&\leq C_\varphi\bigg(\int_{\R^N}u^2(\lambda x,\lambda^2)
  G(x,1)\,dx\bigg)^{\!\!1/2}\bigg(\int_{\R^N}\widetilde V_{m,k}^2(x)
  G(x,1)\,dx\bigg)^{\!\!1/2}\\
  \notag&\qquad+ C_\varphi \bigg(\int_{\R^N}|u(\lambda
  x,\lambda^2)|^{p+1}|G(x,1)|^{\frac{p+1}2}\,dx\bigg)^{\!\!\frac1{p+1}}
  \bigg(\int_{\R^N}|\widetilde
  V_{m,k}(x)|^{p+1}|G(x,1)|^{\frac{p+1}2}\,dx\bigg)^{\!\!\frac1{p+1}}\\
  &\notag\hskip8cm\times \bigg(\int_{\R^N}|u(\lambda
  x,\lambda^2)|^{p+1}\,dx\bigg)^{\!\!\frac{p-1}{p+1}}\\
&\notag\leq C_\varphi\sqrt{H(\lambda^2)}+
C_{\varphi}C_{p+1}\lambda^{-N\frac{p-1}{p+1}}\|u_\lambda(\cdot,1)\|_{\mathcal H}
\|\widetilde
  V_{m,k}\|_{\mathcal H}\bigg(\int_{\R^N}|u(y,\lambda^2)|^{p+1}\,dy
\bigg)^{\!\!\frac{p-1}{p+1}}.
\end{align}
From Corollary \ref{c:pos_per_nonlin} and Lemma \ref{l:Nabove}, it follows that
\begin{multline}\label{eq:92}
\|u_\lambda(\cdot,1)\|_{\mathcal H}^2=\|u(\cdot,\lambda^2)\|_{\mathcal H_{\lambda^2}}^2
\leq  \frac{\lambda^2}{C_1''}\bigg(D(\lambda^2)+
  \frac{N-2}{4\lambda^2}H(\lambda^2)\bigg)\\=
  \frac{H(\lambda^2)}{C_1''}\bigg({\mathcal N}(\lambda^2)+\frac{N-2}{4}\bigg)
\leq\frac{
C_{4}\big({\mathcal N}(T_0)+\frac{N-2}4\big)}{C_1''}
H(\lambda^2),
\end{multline}
while, from Corollary \ref{c:pos_def}, for all $(m,k)\in J_0$,
\begin{equation}\label{eq:95}
\|\widetilde
  V_{m,k}\|_{\mathcal H}\leq
{\rm const\,}\bigg(\gamma+\frac{N-2}4
\bigg).
\end{equation}
From (\ref{eq:94}), (\ref{eq:92}), (\ref{eq:95}),  and Lemma \ref{stimaH},
we deduce that
\begin{align}\label{eq:91}
|\xi_{m,k}(\lambda)|\leq C_7 \lambda^{-2+\frac{N+2-p(N-2)}{p+1}+2\gamma},
\quad\text{for all }\lambda\in(0,\sqrt{T_0})
\end{align}
and for some positive constant $C_7$ depending on
$\|u\|_{L^{\infty}(0,T, L^{p+1}(\R^N))}$, $a$, $N$, $\gamma$, $\varphi$, $T_0$,
$K_1$, $p$,
but independent of $\lambda$ and $(m,k)\in J_0$.
\end{description}
Collecting (\ref{eq:57}) and (\ref{eq:91}), we have that
\begin{align}\label{eq:93}
|\xi_{m,k}(\lambda)|\leq C_8 \lambda^{-2+\tilde\delta+2\gamma},
\quad\text{for all }\lambda\in(0,\sqrt{T_0})
\end{align}
for some $C_8>0$ which is independent of $\lambda$ and $(m,k)\in J_0$ and
$$
\tilde\delta=
\begin{cases}
  \e/2,
  &\text{in case {\bf (I)}},\\[5pt]
  \frac {N+2-p(N-2)}{p+1}, &\text{in case {\bf (II)}}.
\end{cases}
$$
Estimate (\ref{eq:93}) implies that the function $s\mapsto
s^{1-2\gamma}\xi_{m,k}(s)$ belongs to $L^1(0,\sqrt{T_0})$. Therefore,
letting $\bar \lambda\to 0^+$ in (\ref{eq:59}) and using (\ref{eq:58}), we
deduce that,  for all
$\lambda\in(0,\sqrt{T_0})$,
\begin{equation}\label{eq:60}
u_{m,k}(\lambda)=
-2\lambda^{2\gamma}\int_0^\lambda s^{1-2\gamma}\xi_{m,k}(s)\,ds.
\end{equation}
From (\ref{eq:93}) and (\ref{eq:60}), we obtain that, for all
$(m,k)\in J_0$ and $\lambda\in(0,\sqrt{T_0})$,
\begin{equation}\label{eq:61}
|u_{m,k}(\lambda)|\leq \frac{2C_8}{\tilde\delta}\lambda^{2\gamma+\tilde\delta}.
\end{equation}
Let us fix $\sigma\in \big(0,\tilde\delta)$; by Lemma \ref{stimaH},
there exists $K_2(\sigma)$ such that
$$
H(\lambda^2)\geq K_2(\sigma)\lambda^{2(2\gamma+\sigma)}
\quad\text{for }\lambda\in(0,\sqrt{T_0}).
$$
Therefore, in view of (\ref{eq:61}), for all
$(m,k)\in J_0$ and $\lambda\in(0,\sqrt{T_0})$,
$$
\frac{|u_{m,k}(\lambda)|}{\sqrt{H(\lambda^2)}}\leq
\frac{2C_8}{\tilde\delta \sqrt{K_2(\sigma)}}\lambda^{\tilde\delta-\sigma}
$$
and hence
\begin{equation}\label{eq:62}
\frac{u_{m,k}(\lambda)}{\sqrt{H(\lambda^2)}}\to 0\quad
\text{as }\lambda\to 0^+.
\end{equation}
On the other hand, by Lemma \ref{l:blow_up}, for every sequence
$\lambda_n\to 0^+$, there exists a subsequence
$\{\lambda_{n_j}\}_{j\in\N}$ and an eigenfunction $g\in
E_0\setminus\{0\}$ of the operator $L$ associated to $\gamma$ such
that
\begin{equation*}
\frac {u_{\lambda_{n_j}}(x,1)}{\sqrt{H(\lambda_{n_j}^2)}}
\to g\quad \text{in }{\mathcal L}\quad\text{as }j\to+\infty,
\end{equation*}
thus implying, for all $(m,k)\in J_0$,
\begin{equation}\label{eq:63}
\frac{u_{m,k}(\lambda_{n_j})}{\sqrt{H(\lambda_{n_j}^2)}}
=\left(\frac {u_{\lambda_{n_j}}(x,1)}{\sqrt{H(\lambda_{n_j}^2)}},
\widetilde V_{m,k}\right)_{\mathcal L}\to(g,
\widetilde V_{m,k})_{\mathcal L}\quad\text{as }j\to+\infty.
\end{equation}
From (\ref{eq:62}) and (\ref{eq:63}), we deduce that
 $(g,\widetilde V_{m,k})_{\mathcal L}=0$ for all $(m,k)\in J_0$.
Since $g\in E_0$ and $\{\widetilde V_{m,k} : (m,k)\in J_0\}$
is  an orthonormal basis of $E_0$, this implies that $g=0$,
a contradiction.
\end{pf}

We now complete the description of the asymptotics of solutions by
combining Lemmas \ref{l:blow_up} and \ref{limite_pos} and obtaining
some convergence of the blowed-up solution continuously as $\lambda\to
0^+$ and not only along subsequences, thus proving Theorems \ref{asym1}
and \ref{asym2}.

\medskip
\begin{pfn}{Theorems \ref{asym1} and \ref{asym2}}
  Identities (\ref{eq:671}) and (\ref{eq:672}) follow from part (i) of Lemma
  \ref{l:blow_up} and Proposition \ref{p:explicit_spectrum}, 
which imply that
  there exists an eigenvalue
  $\gamma_{m_0,k_0}=m_0-\frac{\alpha_{k_0}}2$ of $L$, $m_0,k_0\in\N$,
  $k_0\geq 1$,  such that $\gamma=\lim_{t\to
    0^+}{\mathcal N}(t)=\gamma_{m_0,k_0}$.  Let $E_0$ be the associated
  eigenspace and $J_0$ the finite set of indices
  $\{(m,k)\in\N\times(\N\setminus\{0\}):\gamma_{m_0,k_0}=m-\frac{\alpha_{k}}2\}$,
  so that $\{\widetilde V_{m,k} : (m,k)\in J_0\}$, with the
  $\widetilde V_{m,k}$'s as in Remark \ref{rem:ortho},  is an
  orthonormal basis of $E_0$. Let $\{\lambda_n\}_{n\in\N}\subset
  (0,+\infty)$ such that $\lim_{n\to+\infty}\lambda_n=0$. Then, from
  part (ii) of Lemma \ref{l:blow_up} and Lemmas \ref{l:limite} and
  \ref{limite_pos}, there exist a subsequence
  $\{\lambda_{n_k}\}_{k\in\N}$ and real numbers
  $\{\beta_{n,j}:(n,j)\in J_0\}$ such that
  $\beta_{n,j}\neq 0$ for some $(n,j)\in J_0$ and, for any $\tau\in(0,1)$,
\begin{equation}\label{eq:74}
\lim_{k\to+\infty}\int_\tau^1
\bigg\|\lambda_{n_k}^{-2\gamma}u(\lambda_{n_k} x,\lambda_{n_k}^2t)
-t^{\gamma}\sum_{(n,j)\in J_0}\beta_{n,j}\widetilde V_{n,j}(x/\sqrt t)
\bigg\|_{{\mathcal H}_t}^2dt=0
\end{equation}
and
\begin{equation}\label{eq:75}
\lim_{k\to+\infty}\sup_{t\in[\tau,1]}
\bigg\|\lambda_{n_k}^{-2\gamma}u(\lambda_{n_k} x,\lambda_{n_k}^2t)
-t^{\gamma}\sum_{(n,j)\in J_0}\beta_{n,j}\widetilde V_{n,j}(x/\sqrt t)
\bigg\|_{{\mathcal L}_t}=0.
\end{equation}
In particular,
\begin{equation}\label{eq:71}
  \lambda_{n_k}^{-2\gamma}u(\lambda_{n_k} x,\lambda_{n_k}^2)
  \mathop{\longrightarrow}\limits_{k\to+\infty}
  \sum_{(n,j)\in J_0}\beta_{n,j}\widetilde V_{n,j}(x)
  \quad\text{in }{\mathcal L}.
\end{equation}
 We now prove
that the $\beta_{n,j}$'s depend neither on the sequence
$\{\lambda_n\}_{n\in\N}$ nor on its subsequence
$\{\lambda_{n_k}\}_{k\in\N}$.
Let us fix $\Lambda\in(0,\sqrt{T_0})$ and define  $u_{m,i}$ and $\xi_{m,i}$
as in (\ref{eq:68}-\ref{eq:69}).  By expanding 
$u_{\lambda}(x,1)=u(\lambda
  x,\lambda^2)\in{\mathcal L}$ in Fourier series
 as in (\ref{eq:70}),
from
(\ref{eq:71}) it follows that, for any $(m,i)\in J_0$,
\begin{equation}\label{eq:72}
\lambda_{n_k}^{-2\gamma}u_{m,i}(\lambda_{n_k})
\to\sum_{(n,j)\in J_0} \beta_{n,j}\int_{\R^N}
\widetilde V_{n,j}(x)\widetilde V_{m,i}(x)G(x,1)\,dx=\beta_{m,i}
\end{equation}
as $k\to+\infty$.
As deduced in the proof of Lemma \ref{limite_pos} (see (\ref{eq:59})),
for any $(m,i)\in J_0$
and $\lambda\in(0,\Lambda)$ there holds
\begin{align}\label{eq:73}
u_{m,i}(\lambda)=
\lambda^{2\gamma}\left(\Lambda^{-2\gamma} u_{m,i}(\Lambda)+
2\int_{\lambda}^{\Lambda}s^{1-2\gamma}\xi_{m,i}(s)\,ds\right).
\end{align}
Furthermore, arguing again as in Lemma \ref{limite_pos} (see
(\ref{eq:93})), $s\mapsto s^{1-2\gamma}\xi_{m,i}(s)$ belongs to
$L^1(0,\sqrt{T_0})$. Hence, combining (\ref{eq:72}) and (\ref{eq:73}),
we obtain, for every $(m,i)\in J_0$,
\begin{align*}
  \beta_{m,i}&=\Lambda^{-2\gamma} u_{m,i}(\Lambda)+
  2\int_{0}^{\Lambda}s^{1-2\gamma}\xi_{m,i}(s)\,ds\\
  &=\Lambda^{-2\gamma}
  \int_{\R^N}u(\Lambda x,\Lambda^2)\widetilde V_{m,i}(x)G(x,1)\,dx\\
  &\hskip3cm+2\int_0^{\Lambda}s^{1-2\gamma} \bigg( \int_{\R^N}f(s
  x,s^2, u(sx,s^2))\widetilde V_{m,i}(x)G(x,1)\,dx \bigg) ds.
\end{align*}
In particular the $\beta_{m,i}$'s depend neither on the sequence
$\{\lambda_n\}_{n\in\N}$ nor on its subsequence
$\{\lambda_{n_k}\}_{k\in\N}$, thus implying that the convergences in
(\ref{eq:74}) and (\ref{eq:75}) actually hold as $\lambda\to 0^+$
and proving the theorems.
\end{pfn}

\noindent 
The strong unique continuation property is a direct consequence 
of Theorems \ref{asym1} and \ref{asym2}.

\begin{pfn}{Corollary \ref{cor:uniq_cont}}
  Let us assume by contradiction that $u\not\equiv 0$ in $\R^N\times
  (0,T)$ and fix $k\in\N$ such that $k>\gamma$, with
  $\gamma=\gamma_{m_0,k_0}$ as in Theorems \ref{asym1} and
  \ref{asym2}.  From assumption (\ref{eq:uniq_cont}), it follows that,
  for a.e. $(x,t)\in\R^N\times (0,1)$,
\begin{equation}\label{eq:97}
\lim_{\lambda\to 0^+}|\lambda^{-2\gamma}t^{-\gamma}u(\lambda x,\lambda^2 t)|=0.
\end{equation}
On the other hand, from Theorems \ref{asym1} and \ref{asym2}, it
follows that there exists $g\in{\mathcal H}\setminus\{0\}$ such that
$g$ is an eigenfunction of the operator $L$ associated to $\gamma$
and, for all $t\in(0,1)$ and a.e. $x\in\R^N$,
$$
\lambda^{-2\gamma}t^{-\gamma}u(\lambda x,\lambda^2 t)
\to g(x/\sqrt t),
$$
which, in view of (\ref{eq:97}), implies $g\equiv 0$, a contradiction.
\end{pfn}

 \textbf{Aknowledgements.} The authors would like to thank
 Prof. Susanna Terracini  for her interest in their work and
for helpful comments and discussions.

\end{document}